\documentclass[12pt, reqno, twoside]{amsart}
\usepackage{amsmath, amsthm, mathrsfs, amscd, amsfonts, amssymb, graphicx, color}
\usepackage[bookmarksnumbered, colorlinks, plainpages]{hyperref}
\hypersetup{colorlinks=true,linkcolor=red, anchorcolor=green, citecolor=cyan, urlcolor=red, filecolor=magenta, pdftoolbar=true}
\usepackage{enumitem}
\usepackage{enumitem}
\newtheorem{theorem}{Theorem}[section]
\newtheorem{lemma}[theorem]{Lemma}

\newtheorem{corollary}[theorem]{Corollary}
\theoremstyle{definition}

\theoremstyle{remark}
\newtheorem{remark}[theorem]{Remark}
\numberwithin{equation}{section}

\begin{document}

\setcounter{page}{1}

\title[Gradient estimates for the fast diffusion equation]
{The nonlinear fast diffusion equation on smooth metric measure spaces: Hamilton-Souplet-Zhang estimates and a Ricci-Perelman super flow}

\author[A. Taheri]{Ali Taheri}

\author[V. Vahidifar]{Vahideh Vahidifar}

\address{School of Mathematical and  Physical Sciences, 
University of Sussex, Falmer, Brighton, United Kingdom.}
\email{\textcolor[rgb]{0.00,0.00,0.84}{a.taheri@sussex.ac.uk}} 
 
\address{School of Mathematical and  Physical Sciences, 
University of Sussex, Falmer, Brighton, United Kingdom.}
\email{\textcolor[rgb]{0.00,0.00,0.84}{v.vahidifar@sussex.ac.uk}}


\begin{abstract}
This article presents new gradient estimates for positive solutions to the nonlinear fast diffusion 
equation on smooth metric measure spaces, involving the $f$-Laplacian. The gradient estimates of interest 
are mainly of Hamilton-Souplet-Zhang or elliptic type and are proved using different set of methods and 
techniques. Various implications notably to parabolic Liouville type results and characterisation of ancient 
solutions are given. The problem is considered in the general setting where the metric and potential evolve 
under a super flow involving the Bakry-\'Emery $m$-Ricci curvature tensor. The remarkable interplay between 
geometry, nonlinearity, and evolution -- and their intricate roles in the estimates and the maximum 
exponent range of fast diffusion -- is at the core of the investigation.   
\end{abstract}

\maketitle
 
{ 
\hypersetup{linkcolor=black}
\tableofcontents
}

\allowdisplaybreaks

\section{Introduction} \label{sec1}

Our principal aim in this paper is to study the nonlinear fast diffusion equation in the framework of smooth metric measure 
spaces (hereafter denoted SMMS for brevity). There are many motivations and reasons for this study. Some, in particular 
include, understanding the analytic and dynamic behaviour of solutions, the probabilistic and stochastic features of 
diffusion on SMMS, the role of geometry and curvature on the one hand and nonlinearity on the other in forming the 
behaviour and features of solutions, and last but not least, functional and geometric inequalities on SMMS and their close 
ties with the fast diffusion equation that has seen a huge revival of interest and progress in recent years 
({\it see}, e.g., \cite{AGS, Bak, BBDGV, BonBon, Bon-Dol, CCL, DaskKe, Hunag-Li, PLu, VC, Wang, Zhang} 
and the references therein for more).

Towards this end, let $(M,g)$ be a complete Riemannian manifold of dimension $n \ge 2$, $d\mu=e^{-f} dv_g$ 
be a weighted measure associated with the metric $g$ and potential $f$, and let $dv_g$ be the usual Riemannain 
volume measure on $M$. The triple $(M, g, d\mu)$ is referred to as a smooth metric measure space. Our prime interest 
in this paper is the nonlinear fast diffusion equation (FDE for brevity) associated with the triple $(M,g,d\mu)$ given by
\begin{align} \label{eq11}
\square_p u(x,t) = \partial_t u (x,t) - \Delta_f u^p (x,t) = \mathscr N(t,x,u(x,t)), \qquad 0<p<1.  
\end{align}

The context in which we consider this nonlinear FDE is one where the metric tensor $g$ and potential 
$f$ are generally allowed to be time dependent. Such problems have moved to the forefront of research 
in geometric analysis since the work of G.~Perelman \cite{Pe02} and the study of forward and backward 
heat equations on manifolds evolving under (super) Ricci-Hamilton or (super) Ricci-Perelman flows. 
Naturally in this setting the weighted measure, the various curvature tensors and the familiar 
(metric dependent) differential operators are all time dependent which makes the analysis 
more complicated and more challenging. Before discussing this aspect of the problem 
further let us take a moment to describe the different components and ingredients in 
equation \eqref{eq11}

The operator $\Delta_f$ in \eqref{eq11} is the $f$-Laplacian (or weighted Laplacian) acting on functions 
$v \in \mathscr{C}^2(M)$ by the prescription 
\begin{equation} \label{f-Lap-definition}
\Delta_f v = {\rm div}_f \nabla v = e^f {\rm div}(e^{-f} \nabla v) = \Delta v - \langle \nabla f, \nabla v\rangle. 
\end{equation} 
It is a symmetric diffusion operator with respect to the invariant weighted measure $d\mu$ and it arises in a 
variety of contexts ranging from probability theory and geometry to quantum field theory and statistical 
mechanics \cite{Bak, Wang}. The $f$-Laplacian $\Delta_f$ is a natural generalisation of the 
Laplace-Beltrami operator to the SMMS context and it coincides with the latter precisely when the potential 
$f$ is spatially constant.

Referring to \eqref{eq11} again, ${\mathscr N}={\mathscr N}(t,x,u)$ is a sufficiently smooth {\it nonlinearity} 
or forcing term depending on the space-time variables $(x, t)$ and the dependent variable $u>0$. 
The form and structural features of this nonlinearity will evidently have physical and analytical significances 
and one of our aims is to understand the way it appears in the estimates and the manner in which it 
influences the solutions.  To the best of our knowledge our results here are the first to consider the 
role of nonlinearity in the FDE and its role in the estimates and maximum exponent range, even for 
the static case. 
\footnote{By the static case we mean the case where the metric tensor $g$ and potential $f$ are 
time {\it independent} which is naturally a very special case.}

As for the geometry and curvature properties of the triple $(M,g,d\mu)$ we have a hierarchy of second order 
symmetric tensor fields on $M$ defined by ({\it see} \cite{Bak})  
\begin{equation} \label{Ricci-m-f-intro}
{\mathscr Ric}^m_f(g) := {\mathscr Ric}(g) + \nabla\nabla f - \frac{\nabla f \otimes \nabla f}{m-n}, \qquad m \ge n \ge 2, 
\end{equation}
called the Bakry-\`Emery $m$-Ricci curvature tensors. Here ${\mathscr Ric}(g)$ is the usual Riemannain Ricci 
curvature tensor of $g$, $\nabla \nabla f={\rm Hess}(f)$ is the Hessian of $f$, and $m \ge n$ is a 
of dimension type constant which is {\bf not necessarily an integer}. In relation to \eqref{Ricci-m-f-intro}, when 
$m=n$, by convention, $f$ can only be a constant, thus giving ${\mathscr Ric}^m_f(g)={\mathscr Ric}(g)$. 
We also allow for $m = \infty$ in which case by formally passing to the limit in \eqref{Ricci-m-f-intro} we set, 
\begin{equation} \label{Ricf def eq} 
{\mathscr Ric}_f(g) = {\mathscr Ric}(g) + \nabla\nabla f.
\end{equation}

The identity relating the $f$-Laplacian $\Delta_f$ to the Bakry-\'Emery Ricci curvature tensor ${\mathscr Ric}_f(g)$ 
is the weighted Bochner-Weitzenb\"ock formula saying that for $v$ of class ${\mathscr C}^3(M)$, 
\begin{equation} \label{Bochner-1}
\frac{1}{2} \Delta_f |\nabla v|^2 - \langle \nabla v, \nabla \Delta_f v \rangle
= |\nabla\nabla v|^2 + {\mathscr Ric}_f (\nabla v, \nabla v).
\end{equation}

As is readily seen \eqref{Bochner-1} is not immediately applicable to the Bakry-\'Emery $m$-Ricci curvature 
tensor \eqref{Ricci-m-f-intro}. However, by adding and subtracting the symmetric second order rank-one tensor 
$[\nabla f \otimes \nabla f]/(m-n)$ to \eqref{Bochner-1} and an application of the Cauchy-Schwarz inequality 
one obtains the {\it inequality} 
\begin{align} \label{Bochner-2}
\frac{1}{2} \Delta_f |\nabla v|^2 - \langle \nabla v, \nabla \Delta_f v \rangle
&\ge \frac{1}{n} (\Delta v)^2 + {\mathscr Ric}_f (\nabla v, \nabla v) \nonumber \\
&\ge \frac{1}{m} (\Delta_f v)^2 + {\mathscr Ric}^m_f (\nabla v, \nabla v).
\end{align}

As an easy inspection reveals, a lower bound on ${\mathscr Ric}_f^m(g)$, in the sense of symmetric tensors, 
is a stronger condition than the same lower bound on ${\mathscr Ric}_f(g)$, in that, we have the one way 
implication 
\begin{equation}
{\mathscr Ric}_f^m(g) = {\mathscr Ric}_f(g) - \frac{\nabla f \otimes \nabla f}{m-n} 
\ge {\mathsf k} g \implies {\mathscr Ric}_f(g) \ge {\mathsf k} g,
\end{equation} 
but not the other way around (here ${\mathsf k} \in {\mathbb R}$, $m \ge n$ are finite constants). 
This observation is an important one throughout the analysis below especially from the point of 
view of comparison with similar results for linear heat-type equations. Also it is the curvature tensor 
${\mathscr Ric}_f^m(g)$ that plays the main role in the estimates and results here.

In this paper we present ourselves with the task of developing gradient estimates of elliptic type, 
also called Hamilton-Souplet-Zhang type, for positive solutions to \eqref{eq11}. These will be 
established using different sets of ideas and will cover different exponent ranges in the fast diffusion 
regime. The remarkable fact that the maximum exponent range in these estimates turns to be independent of 
the evolution of metric and potential as well as the curvature bounds is in sharp contrast to other 
classes of gradient estimates and one interesting outcome of our analysis.

Following the standard theory of degenerate and singular diffusion ({\it see} \cite{BDM, DaskKe, Vaz}), 
in the FDE and SMMS context here we introduce the second order linear (space-time dependent variable coefficient) 
evolution operator  
\begin{align} \label{Lpv-definition-introduction}
\mathscr L = \partial_t - p u^{p-1}(x,t) \Delta_f, 
\end{align}
where $u=u(x,t)$ is a positive (solution) to \eqref{eq11} and $0<p<1$. This operator plays an important role 
throughout the paper, particularly, in that, by suitably transforming $u$, say to $v$ ({\it see} below), it follows 
from $u$ being a positive solution to \eqref{eq11} that $v$ is in turn a positive solution to a related equation 
involving the linear (heat-type) operator $\mathscr L$. See Sections \ref{sec3} and \ref{sec8} for more and 
the discussion surrounding \eqref{SPR-p-substitute-intro} at the end of this section. The transformation 
taking $u$ to $v$ is sometimes called the {\it pressure} transform.

In this paper we prove two sets of gradient estimates of Hamilton-Souplet-Zhang type for positive solutions 
to \eqref{eq11}. The first estimate formulated in Theorem \ref{thm2-FDE} covers the {\bf super-critical range} 
$p_c < p <1$ where $p_c=p_c(m) = 1-2/m$. The second estimate formulated in Theorem \ref{thm-6.1-EQ-FDE} 
covers the range $p_0<p<1$ where $p_0=p_0(m)$ is the larger of $1/2$ and $1-1/\sqrt{m-1}$, i.e., 
$p_0(m)=\max(1/2, 1-1/\sqrt{m-1})$. It is not difficult to see that this latter range (depending on $m \ge 2$) 
covers both the super-critical range and partly the {\bf very fast diffusion} or {\bf sub-critical range} $0<p<p_c$. 
In what follows we refer to the larger of the ranges $p_c(m)<p<1$ and $p_0(m)<p<1$ as the 
{\bf maximum exponent range} (for that particular $m$). Note that $m$ is not necessarily an integer. 
\footnote{The Riemannian case corresponding to $d\mu=dv_g$, ${\mathscr Ric}_f^m(g)={\mathscr Ric}(g)$, 
$\Delta_f=\Delta$ and $m=n$ (with $n \ge 2$ integer) is where these terminologies on exponents are 
borrowed from. Here the critical exponent $p_c$ is also known as the Aronson-B\'enilan exponent 
\cite{AronBen}.}

Indeed, looking more closely at these exponents and noting $1-1/\sqrt{m-1} \le 1-2/m $ (for $m \ge 2$) 
it is easy to describe the maximum exponent range for different values of $m$ (along with the corresponding 
theorem) in the way below: 
\footnote{There is another important exponent $p=1-1/[\sqrt{2m} -1]$ that lies between $1-1/\sqrt{m-1}$ 
and $1-2/m $ and is intimately linked with Theorem \ref{thm2-FDE}. We do not discuss this explicitly here 
but it lurks in the background in the proof of Theorem \ref{thm2-FDE} and formally appears 
in Corollary \ref {first-cor-last-fde1} (for $s=2$) and Remark \ref{q1-q2-remark-fde1} 
and again in Remark \ref{remark-after-second-cor-last-fde1}.}
\begin{itemize}
\item $2 \le m \le 4$: $p_c \le p_0=1/2 \implies$ $p_c<p<1$ (Theorem \ref{thm2-FDE}).
\item $4 \le m \le 5$: $1-1/\sqrt{m-1} \le 1/2 =p_0 \le p_c \implies$  $p_0<p<1$ (Theorem \ref{thm-6.1-EQ-FDE}).
\item $m \ge 5$: $1/2 \le 1-1/\sqrt{m-1} = p_0 \le p_c \implies$ $p_0<p<1$ (Theorem \ref{thm-6.1-EQ-FDE}).
\end{itemize}

A remarkable outcome of our analysis is that although the evolution of the metric $g$ and potential $f$ 
and the analytic properties of the nonlinearity ${\mathscr N}$ directly influence the estimate and bounds, 
however, neither the nonlinearity nor the evolution of the metric and potential, have any influence on the 
maximum exponent range in the estimates.

The evolution operator \eqref{Lpv-definition-introduction} and the corresponding 
transformation from $u$ to $v$ takes different forms in each of the two sets of estimates and corresponding theorems: 
\begin{itemize}
\item (Section \ref{sec2}: Theorems \ref{thm2-FDE} to \ref{ancient-1}, Section \ref{sec6}) Here for $0 \le p_c(m)<p<1$ 
we have $v= p/(1-p) u^{p-1}$ and subsequently  
\begin{equation}
\mathscr L = \partial_t-(1-p ) v \Delta_f.
\end{equation}
\item (Section \ref{sec7}: Theorems \ref{thm-6.1-EQ-FDE} to \ref{ancient-2}, Section \ref{sec11}) Here for 
$1/2 \le p_0(m)<p<1$ we have $v=u^{p-1/2}$ and subsequently 
\begin{equation}
\mathscr L = \partial_t - pv^{2(p-1)/(2p-1)} \Delta_f.
\end{equation}
\end{itemize}

An interesting fact to highlight here is that in the nonlinear diffusion context studied in this paper, 
the metric and potential super flow inequality  
\begin{equation} \label{SPR-p-substitute-intro}
\frac{1}{2} \partial_t g(x,t) + p u^{p-1}(x,t) {\mathscr Ric}_f^m(g)(x,t) \ge - {\mathsf k} g(x,t)
\end{equation}
presents itself as the substitute for the super Ricci-Perelman flow inequality that links to linear heat-type 
equations where $p=1$ (see \cite{Taheri-GE-1, Taheri-GE-2, TVahNA, TVDiffHar}). As will become 
evident later, this formulation is closely connected with the evolution operator ${\mathscr L}$ defined 
in \eqref{Lpv-definition-introduction} and reduces to the usual Ricci-Perelman super flow when $p=1$. 
Again, to make the distinction between the two sets of estimates established here clearer, we point out 
that, the super flow \eqref{SPR-p-substitute-intro} in each of the two cases above takes the form:  
\begin{itemize}
\item For $0\le p_c(m)<p<1$ with transformation $v= p/(1-p) u^{p-1}$ we have   
\begin{equation} \label{SP-1-intro}
\frac{1}{2} \dfrac{\partial g}{\partial t} (x,t) + (1-p) v(x,t) {\mathscr Ric}^m_f(g)(x,t) 
\ge - \mathsf{k}g(x,t), \qquad {\mathsf k} \ge 0.
\end{equation}
\item For $1/2 \le p_0(m)<p<1$ with transformation $v=u^{p-1/2}$ we have  
\begin{equation} \label{SP-2-intro}
\frac{1}{2} \dfrac{\partial g}{\partial t} (x,t) + p v^{2(p-1)/(2p-1)}(x,t) {\mathscr Ric}^m_f(g)(x,t) 
\ge - \mathsf{k}g(x,t), \qquad {\mathsf k} \ge 0. 
\end{equation}
\end{itemize}

As a further comment, we also note that in the main results and estimates in Theorem \ref{thm2-FDE} 
and Theorem \ref{thm-6.1-EQ-FDE} we shall make use of the lower bounds 
\begin{equation} \label{h-k-bounds-intro-discussion}
\partial_t g (x,t) \ge - 2h g(x,t), \qquad  {\mathscr Ric}^m_f(g)(x,t) \ge - k g(x,t), \\ 
\end{equation}
for suitable constants $h$, $k \ge 0$ in a fixed compact space-time cylinder $Q_{R,T}$ with upper base centered 
at the reference point where the estimate is sought. These bounds on the one hand allow for the use of weighted 
(Laplace) comparison theorems and on the other provide control on the time derivative of geodesic distances 
which are two important ingredients in localisation and the proof of gradient estimates. The constants $h$, 
$k \ge 0$ will appear along with other bounds in the geometry-dependent terms in the ultimate formulation 
of the local estimates. In the static case, where the metric and potential do not depend on time 
(i.e., $\partial_t g \equiv 0$ and $\partial_t f \equiv 0$) we can take $h=0$ and our results are 
immediately seen to cover this special case with the usual assumption of Ricci curvature lower bound 
(space only) of the type ${\mathscr Ric}_f^m(g) \ge - {\mathsf k} g$.

\qquad \\
{\bf Plan of the paper.} The main estimates are presented in Theorem \ref{thm2-FDE}  in Section \ref{sec2} and 
Theorem \ref{thm-6.1-EQ-FDE} in Section \ref{sec7} where both estimates are formulated in their local forms. 
These estimates cover different exponent ranges and their proofs rely on different set of techniques. The global 
forms of these estimates follow by imposing appropriate global bounds on the solution, metric, potential as 
well as the Bakry-\'Emery $m$-Ricci curvature tensor and is formulated in Theorem \ref{thm2-FDE-global} and Theorem 
\ref{thm-6.1-EQ-FDE-global} respectively. The static case (time independent metrics $g$ and potentials $f$) which 
is an important case is treated as a by-product of the above estimates and is presented in Theorem \ref{thm2-FDE-static}
and Theorem \ref{thm-6.1-EQ-FDE-static} respectively. As a nice and useful application of these we are able 
to establish parabolic Liouville results which lead to a characterisation of ancient solutions to \eqref{eq11}. 
These appear in turn in Theorem \ref{ancient-1} and Theorem \ref{ancient-2}. As such Sections \ref{sec2} 
and \ref{sec7} contain the main results of the paper whilst the remaining sections are devoted to 
developing the necessary apparatus and tools for establishing the results and proofs. Sections \ref{sec6} 
and \ref{sec11} present another new set of global bounds and estimates by direct consideration of 
the super flow \eqref{SPR-p-substitute-intro} [in the forms \eqref{SP-1-intro} and \eqref{SP-2-intro}] 
and suitable application of maximum principle. Here $M$ is taken to be closed.

\qquad \\
{\bf Notation.} We write $z=z_+ + z_-$ with $z_+=\max(z, 0)$ and $z_-=\min(z, 0)$. 
Fixing a reference point $x_0 \in M$ we denote by $d=d(x,x_0, t)$ the Riemannian distance between $x$ 
and $x_0$ on $M$ with respect to $g=g(t)$. We write $\varrho=\varrho(x,x_0,t)$ for the geodesic radial 
variable measuring distance between $x$ and origin $x_0$ at time $t$. For a fixed space-time 
reference point $(x_0,t_0)$ and $R>0$, $T>0$ we define the space-time cylinder 
\begin{equation}
Q_{R,T} (x_0,t_0) \equiv \{ (x, t) | d(x, x_0, t) \le R, t_0-T \le t \le t_0 \} \subset M \times [t_0-T, t_0].
\end{equation}
When the metric $g$ is time independent, we denote by $\mathscr{B}_\varrho(x_0) \subset M$ the 
geodesic ball of radius $\varrho>0$ centered at $x_0$. It is evident that in this case we have 
\begin{equation}
Q_{R,T} (x_0,t_0) = \mathscr{B}_R(x_0) \times [t_0-T, t_0] \subset M \times [t_0-T, t_0].
\end{equation}
When the choice of the reference point 
is clear from the context we often abbreviate and write $d(x, t)$, $\varrho(x,t)$ or ${\mathscr B}_\varrho$, 
$Q_{R,T}$ respectively. For a function of multiple variables, we denote its partial derivatives with subscripts 
accordingly. For instance for  $\Gamma=\Gamma(x,u)$ we denote its partial derivatives with respect to 
$x=(x_1, \dots, x_n)$ or $u$ by $\Gamma_x$ or $\Gamma_u$.

\section {A Hamilton-Souplet-Zhang estimate in the range $p_c(m)< p < 1$}
\label{sec2}

In this section we present the first set of gradient estimates for positive solutions to \eqref{eq11}. 
These cover the range $p_c< p<1$. Here $p_c=p_c(m) \ge 0$ is the critical exponent defined by 
$p_c=1-2/m$ where $m \ge 2$. 
Denoting by $\beta_1$, $\beta_2$ the roots of the quadratic polynomial 
$\beta^2 - (2-p)/(1-p)\beta +m/2$, we pick $\beta$ 
such that $0<\beta_1<\beta \le 1<\beta_2$, i.e., $\beta \in (\beta_1, \beta_2)$ and $\beta \le 1$. 
[See the paragraph after \eqref{E-eq2.29-FDE} for more on this.] We write  
\begin{align} \label{Sigma-N-definition}
\Sigma(t,x,v)= p [(1-p)v/p]^{\frac{2-p}{1-p}}
\mathscr N \left(t,x,[(1-p)v/p]^{\frac{1}{p-1}}\right), 
\end {align}
for the scaled nonlinearity resulting from the change from $u$ to $v$ (see below).

\begin{theorem} \label{thm2-FDE}
Let $(M, g, d\mu)$ be a complete smooth metric measure space with $d\mu=e^{-f} dv_g$. 
Assume the metric and potential are time dependent, of class $\mathscr{C}^2$ and that for suitable constants 
$k, h \ge 0$ and $m \ge n$ satisfy ${\mathscr Ric}_f^m (g) \ge -(m-1)k g$, $\partial_t g \ge -2hg$ in the compact 
space-time cylinder $Q_{R,T}$ with $R, T >0$. Let $u$ be a positive solution to \eqref{eq11} with $p_c< p < 1$ 
and let $v =p/(1-p) u^{p-1}$ and $M= \sup_{Q_{R,T}} v$. 
Then there exists $C=C(p,\beta,m)>0$ such that for every $(x,t)$ in $Q_{R/2,T}$ with $t>t_0-T$ we have 
\begin{align}\label{eq-2.1}
\frac{|\nabla v|}{v^{\beta/2}} (x,t) 
\le C \left \{ \begin {array}{ll}
\sqrt h M^{(1-\beta)/2} + \left[ \dfrac{k^{1/4}}{\sqrt {R}} + \dfrac{1}{R} + \sqrt k \right] M^{1-\beta/2} 
\\
\\
+ \dfrac{M^{(1-\beta)/2}}{\sqrt {t-t_0+T}} 
+ \sup_{Q_{R, T}} \left\{\left[\dfrac{|\Sigma_x(t,x,v)|}{v^{(3\beta-2)/2}}\right]^{1/3}\right\} 
\\
\\
+ \sup_{Q_{R, T}}\left\{v^{(1-\beta)/2}\left[\dfrac{ \beta \Sigma(t,x,v)}{v}-2 \Sigma_v(t,x,v) \right]_+^{1/2}\right\}
\end{array}
\right\}.
\end{align}
\end{theorem}

\begin{remark} \label{Sigma-N-u-v-remark}
Since $u$ and $v$ are here related to one-another via $u =[(1-p)v/p]^{1/(p-1)}$ we can write \eqref{Sigma-N-definition} 
as $\Sigma(t,x,v) = p u^{p-2} {\mathscr N}(t,x,u)$.  
\end{remark}

Subject to global bounds we have the following global estimate (in space) that results from passing to the limit 
$R \to \infty$ in \eqref{eq-2.1}.

\begin{theorem} \label{thm2-FDE-global}
Let $(M, g, d\mu)$ be a complete smooth metric measure space with $d\mu=e^{-f} dv_g$. 
Assume the metric and potential are time dependent, of class $\mathscr{C}^2$ and that for suitable constants 
$k, h \ge 0$ and $m \ge n$ satisfy ${\mathscr Ric}_f^m (g) \ge -(m-1)k g$, $\partial_t g \ge -2hg$ on 
$M \times [t_0-T, t_0]$. Let $u$ be a positive solution to \eqref{eq11} with $p_c < p < 1$ 
and let $v =p/(1-p) u^{p-1}$ and $M= \sup v$. 
Then there exists $C=C(p,\beta,m)>0$ such that for every $x \in M$ and $t_0-T<t\le t_0$ we have 
\begin{align}\label{eq-2.1-global}
\frac{|\nabla v|}{v^{\beta/2}} (x,t) 
\le C \left \{ \begin {array}{ll}
\sqrt k M^{1-\beta/2} + \left[ \sqrt h + \dfrac{1}{\sqrt {t-t_0+T}} \right] M^{(1-\beta)/2} 
\\
\\
+ \sup_{M \times [t_0-T, t_0]} \left\{\left[\dfrac{|\Sigma_x(t,x,v)|}{v^{(3\beta-2)/2}}\right]^{1/3}\right\} 
\\
\\
+ \sup_{M \times [t_0-T, t_0]} \left\{v^{(1-\beta)/2} \left[\dfrac{ \beta \Sigma(t,x,v)}{v}-2 \Sigma_v(t,x,v) \right]_+^{1/2}\right\}
\end{array}
\right\}.
\end{align}
\end{theorem}

The static case ($\partial_t g \equiv 0$ and $\partial_t f \equiv 0$) is of enough significance to be formulated 
as a separate corollary. Here we describe the local version. The global version follows by passing to the limit $R\to\infty$.

\begin{theorem} \label{thm2-FDE-static}
Let $(M, g, d\mu)$ be a complete smooth metric measure space with $d\mu=e^{-f} dv_g$ and assume 
${\mathscr Ric}_f^m (g) \ge -(m-1)k g$ in ${\mathscr B}_R$ for some $k \ge 0$, $m \ge n$ 
and $R>0$. Let $u$ be a positive solution to \eqref{eq11} with $p_c < p < 1$ 
and let $v =p/(1-p) u^{p-1}$ and $M= \sup_{Q_{R,T}} v$. 
Then there exists $C=C(p,\beta,m)>0$ such that for every $(x,t)$ in $Q_{R/2,T}$ with $t>t_0-T$ we have 
\begin{align}\label{eq-2.1-static}
\frac{|\nabla v|}{v^{\beta/2}} (x,t) 
\le C \left \{ \begin {array}{ll}
\left[ \dfrac{k^{1/4}}{\sqrt {R}} + \dfrac{1}{R} + \sqrt k \right] M^{1-\beta/2} 
\\
\\
+ \dfrac{M^{(1-\beta)/2}}{\sqrt {t-t_0+T}} 
+ \sup_{Q_{R, T}} \left\{\left[\dfrac{|\Sigma_x(t,x,v)|}{v^{(3\beta-2)/2}}\right]^{1/3}\right\} 
\\
\\
+ \sup_{Q_{R, T}}\left\{v^{(1-\beta)/2}\left[\dfrac{ \beta \Sigma(t,x,v)}{v}-2 \Sigma_v(t,x,v) \right]_+^{1/2}\right\}
\end{array}
\right\}.
\end{align}
\end{theorem}

We end this section with an application of the estimates above to parabolic Liouville-type theorems. 
Here by an ancient solution $u=u(x,t)$ to \eqref{eq11} we mean a solution defined on $M$ for 
all negative times, that is, for all $(x,t)$ with $x \in M$ and $-\infty < t <0$.

\begin{theorem} \label{ancient-1}
Let $(M, g, d\mu)$ be a complete smooth metric measure space with $d\mu=e^{-f}dv_g$ and 
${\mathscr Ric}^m_f(g) \ge 0$. Assume $[(2-p)-\beta(1-p)/2] {\mathscr N}(u)/u - {\mathscr N}_u(u) \ge 0$ 
for some $\beta \in (\beta_1, \beta_2)$ and all $u>0$ where $p_c< p<1$. Then any positive ancient solution 
to the nonlinear fast diffusion equation 
\begin{equation} \label{ancient-equation}
\square_p u = \partial_t u - \Delta_f u^p = {\mathscr N}(u(x,t)), 
\end{equation}
satisfying the growth at infinity $1/u(x,t) = o( [\varrho(x) + \sqrt{|t|}]^{2/[(1-p)(2-\beta)]})$ must be spatially 
constant. If, in addition, ${\mathscr N}(u) \ge a$ for some $a>0$ and all $u>0$ then 
\eqref{ancient-equation} admits no such ancient solutions. 
\end{theorem}

\section{Evolution inequalities ${\bf I}$: $\mathscr L = \partial_t -(1-p) v \Delta_f$ acting on $w = |\nabla v|^2/v^{\beta}$} 
\label{sec3}

In this section we derive evolution identities and inequalities for $w = |\nabla v|^2/v^{\beta}$ 
where $v =p/(1-p) u^{p-1}$. The evolution operator is $\mathscr L = \partial_t -(1-p) v \Delta_f$ and $\beta \in {\mathbb R}$ is fixed. 
Note that in the next few lemmas we assume $0<p<1$ but later in the proof of the estimates in Section \ref{sec4} we 
further restrict this range to $1-2/m=p_c<p<1$.

\begin{lemma}\label{LEM-Lem4.2-FDE}
Let $u$ be a positive solution to \eqref{eq11} with $0<p<1$, set $v =p/(1-p) u^{p-1}$ and 
let $\Sigma=\Sigma(t,x,v)$ be as in \eqref{Sigma-N-definition}. Then $v$ satisfies the evolution equation 
\begin{align}\label{EQ-eq4.1-FDE}
\mathscr L [v] = [\partial_t -(1-p) v \Delta_f]v = - |\nabla v|^2 - \Sigma(t,x,v).
\end{align}
\end{lemma}

\begin {proof}
As here $u$ and $v$ are related via $u = [(1-p)v/p]^{1/(p-1)}$ basic calculations give
\begin{align}
\partial_t u &=\partial_t [(1-p)v/p]^{1/(p-1)} = -(1/p) [(1-p)v/p]^{\frac{2-p}{p-1}} \partial_t v, \nonumber \\ 
\Delta_f u^p &= -(1/p) [(1-p)v/p]^{\frac{2-p}{p-1}}[(1-p) v \Delta v - |\nabla v|^2]
+ [(1-p)v/p]^{\frac{1}{p-1}}\langle \nabla f , \nabla v \rangle. 
\end{align}
Now substituting these ingredients into \eqref{eq11} results in 
\begin{align}
\square_p u = \partial_t u - \Delta_f u^p &= 
-(1/p) [(1-p)v/p]^{\frac{2-p}{p-1}} (\partial_t v - (1-p)v \Delta_f v 
+ |\nabla v|^2) \nonumber \\
&= \mathscr N \left(t,x,[(1-p)v/p]^{\frac{1}{p-1}}\right)
\end{align}
which by recalling \eqref{Sigma-N-definition} and rearranging terms gives at once the desired conclusion. 
\end{proof}

\begin{lemma}\label{Lem-4.3-FDE}
Let $(M, g, d\mu)$ be a complete smooth metric measure space with $d\mu=e^{-f} dv_g$. 
Suppose that the metric and potential are time dependent and of class $\mathscr{C}^2$. 
Let $u$ be a positive solution to \eqref{eq11} with $0<p<1$ and $w = |\nabla v|^2/v^{\beta}$
where $v =p/(1-p) u^{p-1}$ and $\beta \in {\mathbb R}$ is a constant. Then $w$ satisfies the evolution equation
\begin{align}\label{EQ-eq-4.10-FDE}
\mathscr L [w]  
= & -[\partial_t g+2(1-p)v{\mathscr Ric}_f^m(g)]\frac{( \nabla v, \nabla v)}{v^\beta}\\
&-\frac{2(1-p)}{v^{\beta}} \left[ v|\nabla\nabla v|^2
+\frac{v \langle \nabla f , \nabla v \rangle^2}{(m-n)}-|\nabla v|^2 \Delta_f v\right]\nonumber\\
&-2[1-\beta(1-p)] \frac{\langle \nabla v, \nabla |\nabla v|^2 \rangle}{v^\beta}
+\beta \frac{ |\nabla v|^2 \Sigma(t,x,v)}{v^{\beta+1}}\nonumber\\
&-\beta [(1-p)(\beta+1)-1]\frac{|\nabla v|^4}{v^{\beta+1}}
-2 \frac{\langle \nabla v, \nabla\Sigma(t,x,v)\rangle}{v^{\beta}}.\nonumber
\end{align} 
\end{lemma}

\begin{proof}
Starting with $w = |\nabla v|^2/v^{\beta}$ and using equation \eqref{EQ-eq4.1-FDE} for $v$ we can write
\begin{align}\label{eq4.3}
\partial_t w =&~ \partial_t \left[\frac{|\nabla v|^2}{v^\beta}\right]
= 2\frac{\langle \nabla v, \nabla \partial_t v \rangle}{v^\beta}
-\frac{[\partial_t g] ( \nabla v, \nabla v)}{v^\beta}
- \beta\frac{v^{\beta-1}|\nabla v|^2 \partial_t v}{v^{2\beta}} \\
=&~ \frac{2\langle \nabla v ,\nabla [(1-p) v \Delta_f v -|\nabla v|^2 
- \Sigma(t,x,v)]\rangle}{v^\beta} -\frac{[\partial_t g] ( \nabla v, \nabla v)}{v^\beta} \nonumber\\
&-\frac{\beta |\nabla v|^2 [(1-p ) v \Delta_f v -|\nabla v|^2 
- \Sigma(t,x,v)]}{v^{\beta+1}}\nonumber\\
=& ~ 2(1-p)\frac{ |\nabla v|^2 \Delta_f v}{v^\beta}
+2(1-p) \frac{\langle \nabla v, \nabla \Delta_f v \rangle}{v^{\beta -1}}
- 2\frac{\langle \nabla v , \nabla |\nabla v|^2 \rangle}{v^\beta} 
+ \beta \frac{ |\nabla v|^4}{v^{\beta+1}}\nonumber\\
&- 2 \frac{\langle \nabla v, \nabla\Sigma(t,x,v)\rangle}{v^\beta}
-\frac{[\partial_t g] ( \nabla v, \nabla v)}{v^\beta}
- \beta (1-p)\frac{|\nabla v|^2 \Delta_f v}{v^\beta} 
+ \beta \frac{|\nabla v|^2 \Sigma(t,x,v)}{v^{\beta+1}}. \nonumber 
\end{align}
Next using $\nabla w = \nabla |\nabla v|^2/v^\beta- \beta |\nabla v|^2 \nabla v/v^{\beta+1}$ 
and $\Delta_f v = \Delta v - \langle \nabla f, \nabla v\rangle$  we have
\begin{align}\label{eq4.8}
\Delta_f w = \frac{\Delta_f |\nabla v|^2}{v^\beta}
- 2 \beta \frac{\langle \nabla v, \nabla |\nabla v|^2 \rangle}{v^{\beta+1}}
- \beta\frac{|\nabla v|^2 \Delta_f v}{v^{\beta+1}}
+\beta(\beta+1) \frac{ |\nabla v|^4}{v^{\beta+2}}.
\end{align} 
Putting \eqref{eq4.3}-\eqref{eq4.8} together, simplifying and rearranging terms then gives
\begin{align}\label{eq4.9}
[\partial_t-(1-p) v \Delta_f] w 
= &-(1-p) \frac{\Delta_f |\nabla v|^2}{v^{\beta -1}} 
-2[1- \beta(1-p)] \frac{\langle \nabla v, \nabla |\nabla v|^2 \rangle}{v^\beta}\nonumber\\
&- \beta[(1-p)(\beta+1)-1]\frac{|\nabla v|^4}{v^{\beta+1}}
+2(1-p) \frac{|\nabla v|^2 \Delta_f v}{v^{\beta}}\nonumber\\
&+2(1-p) \frac{\langle \nabla v, \nabla \Delta_f v \rangle}{v^{\beta -1}}
-\frac{[\partial_t g] ( \nabla v, \nabla v)}{v^\beta} \nonumber\\
&-2 \frac{\langle \nabla v, \nabla \Sigma(t,x,v)\rangle}{v^{\beta}}
+ \beta \frac{ |\nabla v|^2 \Sigma(t,x,v)}{v^{\beta+1}}.
\end{align}

An application of the weighted Bochner-Weitzenb\"ock formula \eqref{Bochner-1} to the sum of the 
first and fifth terms on the right-hand side in \eqref{eq4.9} leads to 
\begin{align} 
\eqref {eq4.9} = &-\frac{2(1-p)}{v^{\beta -1}}\left[ |\nabla\nabla v|^2
+{\mathscr Ric}_f^m (\nabla v, \nabla v)+\frac{\langle \nabla f , \nabla v \rangle^2}{(m-n)}\right] \\
&-2[1- \beta(1-p)] \frac{\langle \nabla v, \nabla |\nabla v|^2 \rangle}{v^\beta}
-\beta [(1-p)(\beta+1)- 1]\frac{|\nabla v|^4}{v^{\beta+1}}\nonumber\\
&-\frac{[\partial_t g] ( \nabla v, \nabla v)}{v^\beta}+2(1-p) \frac{|\nabla v|^2 \Delta_f v}{v^{\beta}}
-2 \frac{\langle \nabla v, \nabla \Sigma(t,x,v) \rangle}{v^{\beta}}
+\beta \frac{ |\nabla v|^2 \Sigma(t,x,v)}{v^{\beta+1}}. \nonumber 
\end{align}
A rearrangement of terms now gives the desired identity. 
\end{proof}

\begin {lemma}\label{Lem4.4-FDE-pre}
Under the assumptions of Lemma $\ref{Lem-4.3-FDE}$, $w$ satisfies the evolution inequality 
\begin{align}\label{EQ-eq-4.2-FDE-in-Lemma}
\mathscr L [w] 
\le &-[\partial_t g + 2(1-p) v {\mathscr Ric}_f^m (g)] \frac{(\nabla v, \nabla v)}{v^\beta} \nonumber \\
&- 2[1- \beta(1-p)] \langle \nabla v, \nabla w \rangle 
+(1-p) \left[\beta^2 -\frac{2-p}{1-p}\beta +\frac{m}{2} \right] v^{\beta-1} w^2\\
& - 2 \frac{\langle \nabla v, \Sigma_x(t,x,v) \rangle}{v^\beta}
+\left[\frac{ \beta \Sigma(t,x,v)}{v}-2 \Sigma_v(t,x,v) \right]w. \nonumber
\end{align}
\end{lemma}

\begin{proof}
Starting from \eqref{EQ-eq-4.10-FDE}, making note of
$\langle \nabla v, \nabla |\nabla v|^2 \rangle/v^\beta = \langle \nabla v, \nabla w \rangle +\beta |\nabla v|^4/v^{\beta+1}$ 
and $\langle \nabla v, \nabla \Sigma(t,x,v) \rangle 
=\langle \nabla v, \Sigma_x(t,x,v) \rangle + \Sigma_v(t,x,v) |\nabla v|^2$
we have after substitution 
\begin{align}\label{EQ-eq-4.18-FDE}
[\partial_t - (1-p) v\Delta_f] w =
& -[\partial_t g+2(1-p)v{\mathscr Ric}_f^m(g)]\frac{( \nabla v, \nabla v)}{v^\beta}\\
&-\frac{2(1-p)}{v^{\beta}} \left[v|\nabla\nabla v|^2
+\frac{v \langle \nabla f , \nabla v \rangle^2}{(m-n)}-|\nabla v|^2 \Delta_f v\right]\nonumber\\
&-2[1-\beta(1-p)] \langle \nabla v, \nabla w \rangle 
+\left[\beta \frac{ \Sigma(t,x,v)}{v}-2 \Sigma_v(t,x,v)\right]\frac{|\nabla v|^2}{v^{\beta}}\nonumber\\
&+[(1-p)\beta^2-(2-p)\beta]\frac{|\nabla v|^4}{v^{\beta+1}}
-2\frac{\langle \nabla v, \Sigma_x(t,x,v) \rangle}{v^\beta}.\nonumber
\end{align}
Next, using the inequality $|\nabla \nabla v|^2 + \langle \nabla f , \nabla v \rangle^2/(m-n) 
\ge (\Delta_f v)^2/m$, it is seen that 
\begin{align} 
&-\frac{2(1-p)}{v^{\beta -1}} \left[ |\nabla\nabla v|^2 + \frac{\langle \nabla f , \nabla v \rangle^2}{m-n}
- \frac{|\nabla v|^2 \Delta_f v}{v}\right] \le -\frac{2(1-p)}{v^{\beta -1}} \left[  \frac{(\Delta_f v)^2}{m}
- \frac{|\nabla v|^2 \Delta_f v}{v}\right] \nonumber \\
& \qquad \le - \frac{2(1-p)}{v^{\beta -1}} \left[\left( \frac{\Delta_f v}{\sqrt m} - \frac{\sqrt m}{2} \frac{|\nabla v|^2}{v} \right)^2
-\frac{m}{4} \frac{|\nabla v|^4}{v^2} \right] \le \frac{m(1-p)}{2} \frac{|\nabla v|^4}{v^{\beta+1}}.
\end{align}
Therefore substituting the above in \eqref{EQ-eq-4.18-FDE},  
making the substitution $w = |\nabla v|^2/v^{\beta}$ and rearranging terms gives the desired conclusion. 
\end{proof}

\begin{lemma}\label{Lem4.4-FDE}
Under the assumptions of Lemma $\ref{Lem4.4-FDE-pre}$, if for some ${\mathsf k} \ge 0$ 
and $m \ge n$ the metric and potential satisfy the super flow   
\begin{align} \label{flow-inequality-SPR-FDE}
\frac{1}{2} \partial_t g + (1-p) v {\mathscr Ric}_f^m (g) \ge - {\mathsf k} g,  
\end{align}
then $w$ satisfies the evolution inequality 
\begin{align}\label{EQ-eq-4.2-FDE}
\mathscr L [w] 
\le &~2 {\mathsf k} w - 2[1- \beta(1-p)] \langle \nabla v, \nabla w \rangle 
- 2 \frac{\langle \nabla v, \Sigma_x(t,x,v) \rangle}{v^\beta} \\
&+(1-p) \left[\beta^2 -\frac{2-p}{1-p}\beta +\frac{m}{2} \right] v^{\beta-1} w^2 
+\left[\frac{ \beta \Sigma(t,x,v)}{v}-2 \Sigma_v(t,x,v) \right]w. \nonumber
\end{align}
\end{lemma}

\begin{proof}
This follows immediately from Lemma \ref{Lem4.4-FDE-pre} by substituting \eqref{flow-inequality-SPR-FDE} 
into \eqref{EQ-eq-4.2-FDE-in-Lemma} and using $w = |\nabla v|^2/v^{\beta}$. 
\end{proof}

\begin{remark}{
Note that under the assumptions ${\mathscr Ric}_f^m(g) \ge -(m-1) kg$ and $\partial_t g \ge - 2h g$ 
(as in Theorem  \ref{thm2-FDE}) we can replace ${\mathsf k}$ in \eqref{EQ-eq-4.2-FDE} with 
$(1-p)(m-1)kv+h$ or subject to $0<v \le M$ (again as in Theorem \ref{thm2-FDE}) with $(1-p)(m-1)kM+h$.
}
\end{remark}

\section{Localisation in space-time cylinders $Q_{R,T}$ and proof of Theorem \ref{thm2-FDE}}
\label{sec4}

As the gradient estimate in Theorem \ref{thm2-FDE} is a local estimate, its proof is naturally based on localisation and requires 
the construction of suitable cut-off functions. To this end, let us fix a base point $x_0 \in M$, a reference time $t_0 \in {\mathbb R}$, 
and $R, T>0$ and then $\tau \in (t_0-T, t_0]$. The following standard lemma grants the existence of a smooth 
function $\bar{\eta}=\bar \eta(\varrho, t)$ of two real variables $\varrho \ge 0$ and $t_0-T \le t \le t_0$ 
respectively satisfying a set of useful bounds and properties for carrying out the {\it cylindrical} localisation 
procedure. The desired cut-off function $\eta$ will then be constructed from $\bar\eta$ via \eqref{cut-off def}.  
\begin{lemma} \label{phi lemma} 
Fix $t_0 \in {\mathbb R}$ and let $R, T>0$. Given $\tau \in (t_0-T, t_0]$ 
there exists a smooth function $\bar{\eta}:[0,\infty) \times [t_0-T, t_0] \to \mathbb{R}$ such that the following 
properties hold:
\begin{enumerate}[label=$(\roman*)$]
\item ${\rm supp} \, \bar{\eta}(\varrho,t) \subset [0,R] \times [t_0-T, t_0]$  and $0 \leq \bar{\eta}(\varrho,t) \leq 1$ in $[0,R] \times [t_0-T, t_0]$,
\item $\bar{\eta}=1$ in $[0,R/2] \times [\tau, t_0]$ and $\partial \bar{\eta}/\partial \varrho =0$ in $[0,R/2] \times [t_0-T, t_0]$, respectively,
\item there exists $c>0$ such that 
\begin{equation}
\frac{|\partial_t \bar{\eta}|}{\sqrt{\bar\eta}} \le \frac{c}{\tau-t_0+T}, 
\end{equation} 
in $[0,\infty)\times[t_0-T,t_0]$ and $\bar{\eta}(\varrho,t_0-T)=0$ for all $\varrho \in [0,\infty)$.
\item $-c_a \bar{\eta}^a/R \le \partial_\varrho \bar{\eta} \le 0$ and $|\partial_{\varrho \varrho} \bar{\eta}| \le c_a \bar{\eta}^a / R^2$ 
hold on $[0, \infty)\times [t_0-T, t_0]$ for every $0<a<1$ and some $c_a>0$.
\end{enumerate}
\end{lemma}

Having the above lemma in place we now move on to introducing a smooth cylindrical 
cut-off function $\eta=\eta(x,t)$ by setting, for $R>0$, $T>0$ and $t_0-T <\tau \le t_0$, 
\begin{equation} \label{cut-off def}
\eta(x,t) = \bar{\eta}(\varrho(x,t), t), 
\end{equation}
for $(x, t) \in M \times [t_0-T, t_0]$. Referring to \eqref{cut-off def} a straightforward calculation gives 
\begin{itemize}
\item $\partial_t \eta = \partial_\varrho \bar\eta \partial_t \varrho + \partial_t \bar\eta$, 
\item $\nabla \eta = \partial_\varrho \bar\eta \nabla \varrho$, and, 
\item $\Delta_f \eta = \partial_{\varrho \varrho} \bar\eta |\nabla \varrho|^2 
+ \partial_\varrho \bar\eta \Delta_f \varrho$. 
\end{itemize}
As a result we have from the above:  
\begin{align} \label{Lpv-eta-cutoff-equation}
\mathscr L [\eta] &= [\partial_t - (1-p)v \Delta_f] \eta \nonumber \\
&= \partial_\varrho \bar\eta \partial_t \varrho + \partial_t \bar\eta 
- (1-p) v [\partial_{\varrho \varrho} \bar\eta |\nabla \varrho|^2 + \partial_\varrho \bar\eta \Delta_f \varrho] \nonumber \\
&= \partial_\varrho \bar\eta \mathscr L [\varrho] + [\partial_t  - (1-p) v |\nabla \varrho|^2 \partial_{\varrho \varrho}] \bar\eta.
\end{align}

With this introduction we now proceed onto completing the proof of Theorem \ref{thm2-FDE}. 
To this end by considering the localised function $\eta w$ we have
\begin{align} \label{eta-w-Lpv-equation}
\mathscr L [\eta w] = w\mathscr L [\eta]
- 2(1-p)v [\langle\nabla \eta ,\nabla (\eta w) \rangle- |\nabla \eta|^2 w]/\eta
+\eta \mathscr L [w].
\end{align}
Using \eqref{EQ-eq-4.2-FDE} in Lemma \ref{Lem4.4-FDE} and  
$\eta \langle \nabla v, \nabla w \rangle = \langle \nabla v, \nabla (\eta w) \rangle - w \langle \nabla v, \nabla \eta \rangle$ 
we can write
\begin{align}\label{EQ-eq-8.18.FDE}
\mathscr L [\eta w] \le &~
w \mathscr L [\eta]
-2(1-p) v \left\langle \frac{\nabla \eta}{\eta} , \nabla (\eta w) \right\rangle
+ 2(1-p)v \frac{|\nabla \eta|^2}{\eta} w \nonumber\\
&+ 2[(1-p) (m-1)k v + h] \eta w
-2 [1-\beta(1-p)] \langle \nabla v , \nabla (\eta w) \rangle\nonumber\\
&+2 [1-\beta(1-p)] w\langle \nabla v , \nabla \eta \rangle
+ (1-p) \left[ \beta ^2 -\frac{2-p}{1-p} \beta +\frac{m}{2} \right] v^{\beta-1} \eta w^2\nonumber\\
&- 2\eta \frac{\langle \nabla v, \Sigma_x(t,x,v) \rangle}{v^\beta} + \left[\frac{ \beta \Sigma(t,x,v)}{v} -2 \Sigma_v(t,x,v)\right] \eta w.
\end{align}

Assume now that the localised function $\eta w$ achieves its maximum over the compact set 
$\{(x,t) | d(x,x_0, t) \le R, t_0-T \le t \le \tau\} \subset M \times [t_0-T, t_0]$ at the point $(x_1, t_1)$. We assume that $(\eta w)(x_1, t_1) >0$ 
as otherwise the conclusion of the theorem holds true with $w(x, \tau) \le 0$ for all $d(x, x_0, \tau) \le R/2$. Since $\eta \equiv 0$ when 
$t=t_0-T$ by part $(iii)$ in Lemma \ref{phi lemma}, it therefore follows from this strict inequality that $t_1>t_0-T$. Hence at the maximum 
point $(x_1,t_1)$  we have the conditions 
\begin{itemize}
\item $\partial_t (\eta w) \ge 0$, 
\item $\nabla(\eta w) =0$, 
\item $\Delta_f(\eta w) \le 0$,
\end{itemize} 
and subsequently due to $v>0$ and $0<p<1$ we also have 
\begin{equation} 
\mathscr L [\eta w] = [\partial_t-(1-p)v \Delta_f] (\eta w)  \ge 0. 
\end{equation} 
Now returning to \eqref{EQ-eq-8.18.FDE}, making use of the above and rearranging terms 
it is seen that at the maximum point $(x_1, t_1)$ we have 
\begin{align}\label{E-eq2.29-FDE}
- (1-p) \left[ \beta ^2 -\frac{2-p}{1-p} \beta 
+\frac{m}{2} \right] v^{\beta-1} \eta w^2 \le
&~ 2[(1-p)(m-1) k v + h] \eta w\nonumber\\
&+ 2 [1-\beta(1-p)] w\langle \nabla v , \nabla \eta \rangle \nonumber\\
&+ 2(1-p)v \frac{|\nabla \eta|^2}{\eta} w
+ w \mathscr L [\eta]\nonumber\\
&+\left[\frac{ \beta \Sigma(t,x,v)}{v}-2 \Sigma_v(t,x,v) \right] \eta w\nonumber\\
&- 2\eta \frac{\langle \nabla v, \Sigma_x(t,x,v) \rangle}{v^\beta}.
\end{align}

We now need to choose $\beta$ such that firstly $\beta^2 -(2-p)/(1-p) \beta +m/2<0$ 
and secondly that $\beta \le 1$ (the reason for this will become clearer shortly). As for the first condition 
we need $[(2-p)/(1-p)]^2>2m$ so that there are two roots 
$\beta_1< \beta_2$ to the quadratic polynomial in $\beta$ and as for the second condition we need $\beta_1<1$. These two conditions can 
be simultaneously satisfied {\it iff} $p \in (1-2/m,1)$. Picking $\beta \in (\beta_1, \beta_2)$ 
where the roots $\beta_1$, $\beta_2$ now satisfy $0<\beta_1<1<\beta_2$, we can set
$-(1-p)[\beta ^2 -(2-p)/(1-p) \beta +m/2] =2/ \gamma$ with $\gamma > 0$.

Returning to \eqref{E-eq2.29-FDE} and multiplying through by $\gamma v^{1-\beta}$ we arrive at the 
following inequality at the maximum point $(x_1, t_1)$  
\begin{align}\label{EQ-eq-8.20-FDE}
2 \eta w^2 \le&~ 
2 \gamma [(1-p)(m-1) k v^{2-\beta} + h v^{1-\beta}]\eta w\nonumber\\
&+ 2 \gamma [1-\beta(1-p)] v^{1-\beta}w\langle \nabla v , \nabla \eta \rangle\nonumber\\
&+ 2 \gamma (1-p) v^{2-\beta} \frac{|\nabla \eta|^2}{\eta} w
+ \gamma v^{1-\beta} w \mathscr L [\eta]\nonumber\\
&+ \gamma v^{1-\beta}\eta w \left[\frac{ \beta \Sigma(t,x,v)}{v}- 2 \Sigma_v(t,x,v) \right]\nonumber\\
&- 2 \gamma v^{1-\beta} \eta \frac{\langle \nabla v, \Sigma_x(t,x,v) \rangle}{v^\beta}. 
\end{align}

Now the task is to provide suitable upper bounds for each of the terms on the right-hand side of \eqref{EQ-eq-8.20-FDE}.
Towards this end we set $M= \sup_{Q_{R,T}} v$. Then starting from the first line, by recalling 
$0 \le \eta \le 1$ and using Young's inequality we can write 
\begin{align} \label{FDE-proof-1-ezafi-label-zadam}
2 \gamma [(1-p)(m-1) k& v^{2-\beta} + h v^{1-\beta}]\eta w \le
2 \gamma [(1-p)(m-1) k v^{2-\beta} + h v^{1-\beta}] \eta^{1/2} w\nonumber\\
&\le \frac{1}{6} \eta w^2 + C \gamma ^2 [(1-p)^2 (m-1)^2 k^2 v^{4-2\beta} + h^2 v^{2-2\beta}] \nonumber\\
&\le \frac{1}{6} \eta w^2 + C \gamma ^2 [(1-p)^2 (m-1)^2 k^2(\sup_{Q_{R, T}} v)^{4-2\beta} 
+ h^2 (\sup_{Q_{R, T}} v)^{2-2\beta}] \nonumber\\
&\le \frac{1}{6} \eta w^2 + C \gamma ^2 [(1-p)^2 (m-1)^2 k^2 M^{4-2\beta} + h^2 M^{2-2\beta}].
\end{align} 
Note that here $\beta>0$ and we do restrict $\beta \le 1$. 
In particular $4- 2\beta \ge 2$ and $2- 2\beta \ge 0$ and so making use of the upper bound $M$ of $v>0$ on $Q_{R,T}$ 
in \eqref {FDE-proof-1-ezafi-label-zadam} is justified. Moving next to the term on the second line of \eqref{EQ-eq-8.20-FDE}, by 
using the Cauchy-Schwarz inequality, recalling the relation $w = |\nabla v|^2/v^{\beta}$ and utilising $(iv)$ in Lemma \ref{phi lemma}, 
we have 
\begin{align}\label{bound1}
2 \gamma [1-\beta(1-p)] v^{1-\beta}w \langle \nabla v , \nabla \eta \rangle 
& \le 2 \gamma |1-\beta(1-p)| v^{1-\beta}w |\nabla v| |\nabla \eta|\nonumber\\
& \le 2 \gamma |1-\beta(1-p)| v^{1-\beta/2} \eta^{3/4} w^{3/2} \frac{|\nabla \eta|}{\eta^{3/4}}\nonumber\\
& \le \frac{1}{6} (\eta ^{3/4} w^{3/2})^{4/3} + C \gamma^4 [1-\beta(1-p)]^4 \left[\frac{|\nabla \eta|}{\eta^{3/4}}\right]^4 v^{4-2\beta}\nonumber\\
& \le \frac{1}{6} \eta w^2 + \frac{C \gamma^4 [1-\beta(1-p)]^4}{R^4} (\sup_{Q_{R, T}} v)^{4-2\beta}\nonumber\\
& \le \frac{1}{6} \eta w^2 + \frac{C \gamma^4 [1-\beta(1-p)]^4}{R^4} M^{4-2\beta}.
\end{align}

In much the same way for the third term in \eqref{EQ-eq-8.20-FDE}, by using Young's inequality and taking advantage of $(iv)$ 
in Lemma \ref{phi lemma} we can write
\begin{align}
2 \gamma (1-p) v^{2-\beta} \frac{|\nabla \eta|^2}{\eta} w 
&  \le \frac{1}{6} \eta w^2 +C \gamma^2 (1-p)^2 \left[\frac{|\nabla \eta|}{\eta ^{3/4}} \right]^4 v^{4-2\beta}\nonumber\\
& \le \frac{1}{6} \eta w^2 + \frac{C \gamma^2 (1-p)^2}{R^4}(\sup_{Q_{R, T}} v)^{4-2\beta}\nonumber\\
& \le \frac{1}{6} \eta w^2 + \frac{C \gamma^2 (1-p)^2}{R^4} M^{4-2\beta}.
\end{align}

We now come to the term 
$\gamma v^{1-\beta} w \mathscr L [\eta] = \gamma v^{1-\beta} w[\partial_t - (1-p)v \Delta_f] \eta$ which has a purely 
geometric nature. We make use of the curvature bound ${\mathscr Ric}_f^m(g) \ge -(m-1)k g$ and the metric time-derivative 
bound $\partial_t g \ge -2hg$ as given in the theorem to give an upper bound. In fact, to handle this term more efficiently we consider the 
space and time derivatives separately. 
\begin{itemize}
\item The term ${\bf I} = \gamma (1-p) v^{2-\beta} w (-\Delta_f \eta)$: 
By virtue of ${\mathscr Ric}_f^m(g) \geq - (m-1)k g$ and the weighted Laplace comparison theorem we have
$\Delta_f \varrho \le(m-1) \sqrt {k}\coth (\sqrt{k} \varrho)$
and so from \eqref{cut-off def} and the calculations leading to \eqref{Lpv-eta-cutoff-equation} we have  
\begin{align}\label{eq2.26}
\Delta_f \eta = \bar\eta_{\varrho \varrho}|\nabla \varrho|^2 + \bar\eta_\varrho \Delta_f \varrho
\ge \bar\eta_{\varrho \varrho} + \bar\eta_\varrho (m-1)\sqrt{k} \coth (\sqrt{k} \varrho).
\end{align} 
Now using the bound $\sqrt{k} \coth (\sqrt{k} \varrho)\le \sqrt{k} \coth(\sqrt {k} R/2) \le (2+\sqrt {k} R)/R $
for $0 \le \varrho \le R/2$ and noting that $ \bar\eta_\varrho = 0$ for  $0 \le \varrho \le R/2$, this gives
\begin{align}\label{eq2.27}
-\Delta_f \eta& \le - [ \bar\eta_{\varrho \varrho} + \bar\eta_\varrho (m-1)\sqrt {k} \coth(\sqrt {k} R/2)]\nonumber\\
&\le |\bar\eta_{\varrho \varrho}| + (m-1)\left(\frac{2}{R} +\sqrt {k}\right) |\bar\eta_\varrho|.
\end{align} 
Therefore substituting the above bound in the expression for ${\bf I}$ gives 
\begin{align}
{\bf I} &=\gamma (1-p) v^{2-\beta} w (-\Delta_f \eta) \nonumber \\
& \le\gamma (1-p) v^{2-\beta}  \sqrt \eta w 
 \left[\frac{|\bar\eta_{\varrho \varrho}|}{\sqrt \eta}
+(m-1)\left(\frac{2}{R} +\sqrt {k} \right) \frac{|\bar\eta_\varrho|}{\sqrt \eta} \right], 
\end{align}
and so using Young's inequality and $(iv)$ in Lemma \ref{phi lemma} it follows that 
\begin{align}\label{eq2.38}
{\bf I} &\le \frac{\eta w^2}{12} + C\gamma^2 (1-p)^2
\left[\left(\frac{|\bar\eta_{\varrho \varrho}|}{\sqrt \eta}\right)^2
+(m-1)^2\left(\frac{1}{R^2} +k \right) \left(\frac{|\bar\eta_\varrho|}{\sqrt \eta}\right)^2 \right] v^{4-2\beta}\nonumber\\
& \le \frac{\eta w^2}{12} + C \gamma^2(1-p)^2(m -1)^2
\left[\frac{1+ k R^2}{R^4}\right] (\sup_{Q_{R, T}} v)^{4-2\beta}. 
\end{align}

\item The term ${\bf II} = \gamma v^{1-\beta}w \partial_t \eta$:
Here we use the bound $\partial_t \varrho (x,t) \ge - h R$ in $Q_{R,T}$ and 
Lemma \ref{phi lemma}. 
Hence noting $\partial_t \eta = \bar{\eta}_t + \bar \eta_\varrho \partial_t \varrho$ 
and $(iii)$ in Lemma \ref{phi lemma} we have 
\begin{align}\label{eq2.30}
\partial_t \eta &= \bar{\eta}_t + \bar \eta_\varrho \partial_t \varrho 
\le |\bar \eta_t| - h R \bar \eta_\varrho 
= |\bar \eta_t| + h R |\bar \eta_\varrho| \nonumber \\
&\le C \left[ \frac{1}{\tau-t_0+T} + h \right] \sqrt{\eta}.
\end{align} 
Hence substituting in the expression for ${\bf II}$ and using Young's inequality we can write 
\begin{align}
{\bf II} &= \gamma v^{1-\beta}w \partial_t \eta = \gamma \sqrt \eta v^{1-\beta} w \frac{\partial_t \eta}{\sqrt \eta} 
\le C \gamma \sqrt \eta v^{1-\beta} w \left[\frac{1}{\tau-t_0+T} + h \right]\nonumber\\
&\le \frac{1}{12} \eta w^2 + C \gamma ^2 \left[\frac{1}{(\tau-t_0+T)^2} + h^2\right](\sup_{Q_{R, T}} v)^{2-2\beta}.
\end{align}
\end{itemize}
Therefore by putting together the two conclusions above for ${\bf I}$ and ${\bf II}$ we arrive at the bound
\begin{align}
\gamma v^{1-\beta} w \mathscr L [\eta] =&~\overbrace{\gamma (1-p) v^{2-\beta} w (-\Delta_f \eta)}^{{\bf I}} 
+ \overbrace{\gamma v^{1-\beta}w \partial_t \eta}^{{\bf II}} \nonumber \\
\le&~ \frac{1}{6} \eta w^2
+ C \gamma ^2 \left[\frac{1}{(\tau-t_0+T)^2} + h^2\right](\sup_{Q_{R, T}} v)^{2-2\beta}\nonumber\\
&+ C \gamma^2(1-p)^2 (m-1)^2\left[\frac{1+ k R^2}{R^4}\right] (\sup_{Q_{R, T}} v)^{4-2\beta}\nonumber\\
\le&~ \frac{1}{6} \eta w^2 
+ C \gamma ^2 \left\{ \left[\frac{1}{(\tau-t_0+T)^2} + h^2\right] M^{2-2\beta} \right. \nonumber \\
&+ \left. (1-p)^2 (m-1)^2\left[\frac{1+ k R^2}{R^4}\right] M^{4-2\beta} \right\}.
\end{align}

Finally for the two remaining $\Sigma(t,x,v)$ terms, upon recalling $0\le \eta \le 1$ and utilising Young's inequality we have
\begin{align}\label{eq-2.33}
\gamma v^{1-\beta}\eta w & \bigg[\frac{\beta \Sigma(t,x,v)}{v}-2 \Sigma_v(t,x,v) \bigg]
\le \gamma v^{1-\beta}\sqrt \eta w \bigg[\frac{\beta \Sigma(t,x,v)}{v}-2 \Sigma_v(t,x,v) \bigg]_+ \nonumber \\
& \le \frac{1}{6} \eta w^2 
+ C \gamma^2 \sup_{Q_{R, T}} 
\left\{ v^{2(1-\beta)} \left[\frac{ \beta \Sigma(t,x,v)}{v}- 2 \Sigma_v(t,x,v) \right]_+^2\right\},  
\end{align}
and likewise, by using the Cauchy-Schwarz and Young inequalities in turn we can write
\begin{align}
-2 \gamma v^{1-\beta} \eta \frac{\langle \nabla v, \Sigma_x(t,x,v) \rangle}{v^\beta} 
\le &~2 \gamma v^{1-\beta} \eta |\nabla v| \frac{|\Sigma_x(t,x,v)|}{v^\beta} 
\le 2 \gamma \eta^{1/4} w^{1/2} v^{1- \beta/2} \frac{|\Sigma_x(t,x,v)|}{v^\beta} \nonumber\\
\le &~ \frac{1}{6} \eta w^2+C \gamma^{4/3} \sup_{Q_{R, T}} 
\left\{v^{(4/3)[1-\beta/2]} \left[\frac{|\Sigma_x(t,x,v)|}{v^\beta}\right]^{4/3}\right\}. 
\end{align}

Substituting all the above bounds in \eqref{EQ-eq-8.20-FDE}, rearranging terms and making note 
of $0 < \beta \le 1$ it follows that 
\begin{align}
\eta w^2 
\le &~C(\gamma) \left \{ \left( 
\begin {array}{ll} 
(1-p)^2 (m-1)^2\left[\dfrac{1+ k R^2}{R^4}\right]
+\dfrac{(1-p)^2}{R^4}
\\
\\
+\dfrac{ [1-\beta(1-p)]^4}{R^4}
+(1-p)^2 (m-1)^2 k^2  
\end{array}
\right) \right. M^{4-2\beta}
\nonumber \\ 
& \qquad + 
\left. \left( 
\begin {array}{ll} 
\left[ \dfrac{1}{(\tau-t_0+T)^2}+ 2h^2 \right] M^{2-2\beta} 
+ \sup_{Q_{R, T}} \left\{\left[\dfrac{|\Sigma_x(t,x,v)|}{v^{(3\beta-2)/2}}\right]^{4/3}\right\}
\\
\\
+\sup_{Q_{R, T}} \left\{v^{2(1-\beta)}\left[\dfrac{\beta \Sigma(t,x,v)}{v}- 2 \Sigma_v(t,x,v) \right]_+^2\right\}
\end{array}
\right) \right\}.
\end{align}

Recalling the maximality of $\eta w$ at $(x_1, t_1)$ along with $\eta \equiv 1$ when $d(x, x_0,t) \le R/2$ and $\tau \le t \le t_0$, 
it follows that $w^2 (x, \tau) = (\eta^2 w^2)(x, \tau) \le (\eta^2 w^2)(x_1,t_1) \leq (\eta w^2)(x_1,t_1)$. 
Hence noting $w = |\nabla v|^2/v^{\beta}$ and 
absorbing $\beta$, $p$ and $m$ into $C(\gamma)$ gives
\begin{align}
\frac{|\nabla v|}{v^{\beta/2}} (x, \tau)
\le C(p, \beta, m) \left \{ \begin {array}{ll}
\left[ \dfrac{k^{1/4}}{\sqrt {R}}+\dfrac{1}{R}+\sqrt k \right] M^{1-\beta/2} + \sqrt h M^{(1-\beta)/2}
\\
\\
+ \dfrac{M^{(1-\beta)/2}}{\sqrt {\tau-t_0+T}} 
+ \sup_{Q_{R, T}} \left\{\left[\dfrac{|\Sigma_x(t,x,v)|}{v^{(3\beta-2)/2}}\right]^{1/3}\right\}
\\
\\
+ \sup_{Q_{R, T}} \left\{v^{(1-\beta)/2} \left[\dfrac{\beta \Sigma(t,x,v)}{v}- 2 \Sigma_v(t,x,v) \right]_+^{1/2}\right\}
\end{array}
\right\}.
\end{align}
The arbitrariness of $\tau$ in the interval $(t_0-T, t_0]$ now gives the desired conclusion.
\hfill $\square$

\section{Ancient solutions to $\partial_t u - \Delta_f u^p = {\mathscr N}(u)$ ${\bf I}$: Proof of Theorem \ref{ancient-1}}
\label{sec5}

Let us now attend to the proof of Theorem \ref{ancient-1}. Towards this end fix a space-time point $(x_0,t_0)$ and set 
$t=t_0$, $R>0$ and $T=R^2$. It then follows from the local estimate \eqref{eq-2.1-static} in Theorem \ref{thm2-FDE-static} 
(with $k=0$) that 
\begin{align} \label{vx0t0-ancient-1}
\frac{|\nabla v|}{v^{\beta/2}}(x_0,t_0) 
&\le C \left \{ \begin {array}{ll}
\dfrac{M^{1-\beta/2}}{R} + \dfrac{M^{(1-\beta)/2}}{\sqrt T} 
+ \sup_{Q_{R, T}} \left\{\left[\dfrac{|\Sigma_x(t,x,v)|}{v^{(3\beta-2)/2}}\right]^{1/3}\right\} 
\\
\\
+ \sup_{Q_{R, T}}\left\{v^{(1-\beta)/2} \left[\dfrac{ \beta \Sigma(t,x,v)}{v} - 2 \Sigma_v(t,x,v) \right]_+^{1/2}\right\}
\end{array}
\right\}.
\end{align}

By Remark \ref{Sigma-N-u-v-remark} and \eqref{ancient-equation} we have $\Sigma(t,x,v) = p u^{p-2} {\mathscr N}(u)$ 
with $u=[(1-p)v/p]^{1/(p-1)}$. Hence $\Sigma$ does not depend explicitly on $t$ and $x$ and so in particular 
$\Sigma_x \equiv 0$. Moreover by writing $\Sigma_v=\Sigma_u \partial_v u$ where 
$\Sigma_u = pu^{p-2}[(p-2){\mathscr N}/u+{\mathscr N}_u]$ and $\partial_v u = -u^{2-p}/p$ we have  
$\Sigma_v=(2-p){\mathscr N}/u-{\mathscr N}_u$. There using the assumptions in the theorem   
\begin{align}
\beta \frac{\Sigma(v)}{v} - 2 \Sigma_v (v) 
&= \beta (1-p) \frac{{\mathscr N}(u)}{u} - 2 \left[ (2-p) \frac{{\mathscr N}(u)}{u} - {\mathscr N}_u(u) \right] \nonumber \\ 
&= [\beta(1-p)-2(2-p)] \frac{{\mathscr N}(u)}{u} + 2 {\mathscr N}_u(u) \le 0,  
\end{align}
and so $[\beta \Sigma/v-2\Sigma_v]_+\equiv 0$. Now from the stated growth assumption in the theorem 
$1/u(x,t) = o([\varrho(x) + \sqrt{|t|}]^{2/[(1-p)(2-\beta)]})$ we obtain  
\begin{equation} \label{M-growth-ancient-1}
v = [p/(1-p)] u^{p-1} = o([\varrho(x) + \sqrt{|t|}]^{2/(2-\beta)}), 
\end{equation} 
and as a result $M = \sup_{Q_{R,T}} v = o(R^{2/(2-\beta)})$. Therefore from the above and \eqref {vx0t0-ancient-1} we conclude that  
\begin{align} 
\frac{|\nabla v|}{v^{\beta/2}}(x_0,t_0) \le C \left[ \dfrac{M^{1-\beta/2}}{R} + \dfrac{M^{(1-\beta)/2}}{\sqrt T} \right] 
\le \frac{o(R)}{R} + \frac{o(R^{(1-\beta)/(2-\beta)})}{R}.
\end{align}
Now passing to the limit $R \nearrow \infty$ it follows that $|\nabla v|(x_0,t_0)=0$. The arbitrariness of $(x_0,t_0)$ implies 
$|\nabla v| \equiv 0$ and so $v$ and subsequently $u$ are spatially constant. Hence we have $u=u(t)$. From the equation 
satisfied by $u$ it then follows that $du/dt = {\mathscr N}(u)$. Assuming now that ${\mathscr N}(u) \ge a >0$ for all 
$u>0$ it then follows by integrating the ODE that $u(t) \le u(0) + at$ for all $t<0$. This however clashes with 
$u(t)>0$ as $t \searrow -\infty$ and so the conclusion is reached. \hfill $\square$

\section{The Ricci-Perelman super flow ${\bf I}$: $\partial_t g + 2(1-p) v {\mathscr Ric}^m_f(g) \ge -2\mathsf{k}g$} 
\label{sec6}

In this section we obtain further global gradient bounds on positive smooth solutions 
to \eqref{eq11} in the {\it super-critical} range when $M$ is closed and the metric and 
potential evolve under the super flow [{\it cf.} \eqref{SPR-p-substitute-intro} and \eqref{SP-1-intro}]
\begin{align} \label{eq11c-f} 
\frac{1}{2} \dfrac{\partial g}{\partial t} (x,t) + (1-p) v(x,t) {\mathscr Ric}^m_f(g)(x,t) 
\ge - \mathsf{k}g(x,t), \qquad {\mathsf k} \ge 0.
\end{align}

\begin{lemma} \label{evolution Rpq-f}
Let u be a positive solution to \eqref{eq11} and set $v =p/(1-p) u^{p-1}$ where $0<p<1$. 
For $s \ge 2$, $q \in {\mathbb R}$, $\zeta=\zeta(t)$ non-negative and of class $\mathscr{C}^1[0, \infty)$ let
\begin{equation} \label{Xpq def-f}
{\mathsf H}^{s,q}_\zeta[v] = \zeta(t) \frac{|\nabla v|^s}{v^q} + \Gamma(v).
\end{equation} 
Here $\Gamma=\Gamma(v)$ is a function of class $\mathscr{C}^2(0, \infty)$. 
Then $\mathsf H = {\mathsf H}^{s, q}_\zeta [v]$ satisfies the evolution identity,  
\begin{align} \label{R pq phi equality-f}
{\mathscr L} ({\mathsf H}^{s,q}_\zeta[v]) 
=&~ \zeta'(t) \frac{|\nabla v|^s}{v^q} - s \zeta (t) \frac{|\nabla v|^{s-2}}{v^q} 
\left[ \frac{1}{2} \partial_t g + (1-p) v {\mathscr Ric}_f^m(g) \right] ( \nabla v, \nabla v) \\
&+s(1-p)\zeta(t) \frac{|\nabla v|^{s-2}}{v^q} \left[|\nabla v|^2 \Delta_f v - v |\nabla \nabla v|^2
- v \frac{\langle \nabla f , \nabla v \rangle^2}{(m-n)} \right]\nonumber\\
&+ s [q(1-p)-1] \zeta(t) \langle \nabla |\nabla v|^2 , \nabla v\rangle \frac{|\nabla v|^{s-2}}{v^{q}}
-q[q(1-p) -p] \zeta(t) \frac{|\nabla v|^{s+2}}{v^{q+1}}\nonumber\\
&-\frac{s(1-p)}{2v^{q-1}} \zeta(t) \langle \nabla |\nabla v|^{s-2} , \nabla |\nabla v|^2 \rangle
- s \zeta (t) \frac{|\nabla v|^{s-2}}{v^q}\langle \nabla v, \nabla \Sigma(t,x,v) \rangle\nonumber\\
& + q \zeta (t)\frac{|\nabla v|^s}{v^{q+1}} \Sigma(t,x,v) - \Gamma'(v) \Sigma(t,x,v)
- [\Gamma'(v) + (1-p)v\Gamma''(v)]  |\nabla v|^2. \nonumber 
\end{align}
\end{lemma}

\begin{proof}
Referring to \eqref{Xpq def-f} it is seen in view of linearity that
\begin{align}\label{eq-10.5}
{\mathscr L} ({\mathsf H}^{s,q}_\zeta[v]) = 
\zeta (t) {\mathscr L} \left[\frac{|\nabla v|^s}{v^q}\right] 
+ \zeta'(t) \frac{|\nabla v|^s}{v^q} 
+ {\mathscr L} [\Gamma(v)].
\end{align}

We need to evaluate the expressions in the first and third terms on the right-hand
side respectively. Focusing on the first term, it is easily seen that,
\begin{align}\label{10.6.n-f}
{\mathscr L} \left[\frac{|\nabla v|^s}{v^q} \right] 
=&~ \frac{1}{v^q} {\mathscr L} [|\nabla v|^s] 
+ 2 (1-p) v \left \langle \nabla \left(\frac{|\nabla v|^s}{v^q}\right), \frac{\nabla v^q}{v^q} \right \rangle 
- \frac{|\nabla v|^s}{v^{2q}}  {\mathscr L} [v^q].
\end{align}
Now regarding the first term on the right-hand side of \eqref{10.6.n-f}, a direct calculation gives
\begin{align}\label{EQ-eq10.7.n-f}
{\mathscr L} [|\nabla v|^s] 
=&~\frac{s}{2} |\nabla v|^{s-2} \left[ \partial_t -(1-p)v \Delta_f\right] |\nabla v|^2
- \frac{s}{2}(1-p)v \langle \nabla |\nabla v|^{s-2} , \nabla |\nabla v|^2 \rangle \nonumber\\
=&~ \frac{s}{2} |\nabla v|^{s-2} {\mathscr L} [|\nabla v|^2] 
- \frac{s}{2}(1-p)v \langle \nabla |\nabla v|^{s-2} , \nabla |\nabla v|^2 \rangle,
\end{align}
whilst using Lemma \ref {LEM-Lem4.2-FDE} and the weighted Bochner-Weitzenb\"ock formula \eqref{Bochner-1} we have
\begin{align}
{\mathscr L} [|\nabla v|^2]
= & -[\partial_t g] ( \nabla v, \nabla v) +2 (1-p) |\nabla v|^2 \Delta_f v - 2\langle \nabla v, \nabla |\nabla v|^2 \rangle\nonumber\\
&-2 \langle \nabla v, \nabla \Sigma(t,x,v) \rangle-2(1-p) v |\nabla \nabla v|^2\nonumber\\
& -2(1-p)v {\mathscr Ric}_f^m (\nabla v, \nabla v) - 2 (1-p) v \langle \nabla f , \nabla v \rangle^2/(m-n).
\end{align}
Therefore by combining the two identities above it follows that
\begin{align}\label{EQ-eq-2.8}
{\mathscr L} [|\nabla v|^s] 
=&-\frac{s}{2} |\nabla v|^{s-2} [\partial_t g] ( \nabla v, \nabla v) +  s (1-p) |\nabla v|^s \Delta_f v
- s  |\nabla v|^{s-2}\langle \nabla v, \nabla |\nabla v|^2 \rangle\nonumber\\
&- s |\nabla v|^{s-2}\langle \nabla v, \nabla \Sigma(t,x,v) \rangle- s (1-p) v |\nabla v|^{s-2} |\nabla \nabla v|^2\nonumber\\
& - s (1-p)v |\nabla v|^{s-2} {\mathscr Ric}_f^m (\nabla v, \nabla v)
-s (1-p)v |\nabla v|^{s-2} \langle \nabla f , \nabla v \rangle^2/(m-n)\nonumber\\
&-\frac{s}{2}(1-p)v \langle \nabla |\nabla v|^{s-2} , \nabla |\nabla v|^2 \rangle.
\end{align}
Next, for the second term on the right-hand side of \eqref{10.6.n-f} we can write
\begin{align}\label{10.10.n-f}
\left \langle \nabla \left(\frac{|\nabla v|^s}{v^q}\right), \frac{\nabla v^q}{v^q} \right \rangle
= \frac{sq}{2} \frac{|\nabla v|^{s-2}}{v^{q+1}}\langle \nabla |\nabla v|^2 , \nabla v\rangle  - q^2 \frac{|\nabla v|^{s+2}}{v^{q+2}},
\end{align}
and in a similar way for the third term on the right-hand side of \eqref{10.6.n-f} we have 
\begin{align}\label{10.12.n-f}
{\mathscr L} [v^q] 
&= q v^{q-1} [\partial_t - (1-p) v \Delta_f] v - q(1-p)(q-1) v^{q-1} |\nabla v|^2\nonumber\\
&=-q v^{q-1} |\nabla v|^2- q v^{q-1}\Sigma(t,x,v) - q (1-p)(q-1) v^{q-1} |\nabla v|^2 \nonumber \\
&=-q v^{q-1}\Sigma(t,x,v) - q [q(1-p)+p] v^{q-1} |\nabla v|^2. 
\end{align}
Substituting \eqref{EQ-eq-2.8}, \eqref{10.10.n-f} and \eqref{10.12.n-f} back into \eqref{10.6.n-f} 
and rearranging terms leads to
\begin{align}\label{eq-10.14-f}
{\mathscr L} \left[ \frac{|\nabla v|^s}{v^q} \right] 
=& -s\frac{|\nabla v|^{s-2}}{v^q} \left[\frac{1}{2}\partial_t g+(1-p) v {\mathscr Ric}_f^m(g)\right] ( \nabla v, \nabla v) \nonumber\\
&+ s (1-p) \frac{|\nabla v|^{s-2}}{v^q} \left[|\nabla v|^2 \Delta_f v - v |\nabla \nabla v|^2
- v \frac{\langle \nabla f , \nabla v \rangle^2}{(m-n)} \right]\nonumber\\
&- s[1-q(1-p)] \langle \nabla |\nabla v|^2 , \nabla v\rangle \frac{|\nabla v|^{s-2}}{v^{q}}
- q[q(1-p) -p] \frac{|\nabla v|^{s+2}}{v^{q+1}}\nonumber\\
&-\frac{s(1-p)}{2v^{q-1}} \langle \nabla |\nabla v|^{s-2} , \nabla |\nabla v|^2 \rangle
- s \frac{|\nabla v|^{s-2}}{v^q}\langle \nabla v, \nabla \Sigma(t,x,v) \rangle \nonumber\\
& + \frac{q|\nabla v|^s}{v^{q+1}} \Sigma(t,x,v).
\end{align}
 
This completes the calculation of the first term on the right-hand side of \eqref{eq-10.5}.
Since for the third term of the same equation we can write
 \begin{align}\label{eq-10.15-f}
{\mathscr L} (\Gamma(v)) 
&= \Gamma'(v) [-|\nabla v|^2 - \Sigma(t,x,v)] - (1-p)v\Gamma''(v)  |\nabla v|^2,
\end{align}
the conclusion follows upon substituting \eqref{eq-10.14-f} and \eqref{eq-10.15-f} back into \eqref{eq-10.5}.
\end{proof}

\begin{lemma} \label{ERAs-PandS-group}
Under the assumptions of Lemma $\ref{evolution Rpq-f}$, if the metric and potential evolve under the 
super flow \eqref{eq11c-f} then $\mathsf H = {\mathsf H}^{s, q}_\zeta [v]$ satisfies the evolution inequality 
\begin{align} \label{EQ-eq-4.4-2.n-2-fde1}
{\mathscr L} ({\mathsf H}^{s,q}_\zeta[v]) - 2 &[q(1-p)-1] \langle \nabla v, \nabla {\mathsf H}^{s,q}_\zeta[v]\rangle \nonumber\\
\le&~ \frac{|\nabla v|^s}{v^q} \left\{ \zeta'(t) + s {\mathsf k} \zeta(t) + \zeta(t) \left[ q \frac{\Sigma(t,x,v)}{v} - s \Sigma_v(t,x,v)] \right] \right\} \nonumber\\
&+ \zeta(t) (1-p) \left[q^2 - \frac{2-p}{1-p} q+\frac{sm}{4} \right] \frac{|\nabla v|^{s+2}}{v^{q+1}}\nonumber\\
&- s \zeta (t) \frac{|\nabla v|^{s-2}}{v^q}\langle \nabla v, \Sigma_x (t,x,v) \rangle
+ \Gamma'(v) \Sigma(t,x,v) \nonumber\\
& - (1-p) \left[ \frac{2q(1-p)-1}{1-p} \Gamma'(v) + v\Gamma''(v) \right]  |\nabla v|^2. 
\end{align}
\end{lemma}

\begin{proof}
Referring to the second line in the right-hand side of \eqref{R pq phi equality-f} and making note of the inequality 
$|\nabla \nabla v|^2 + \langle \nabla f , \nabla v \rangle^2/(m-n) \ge (\Delta_f v)^2/m$, we can write 
\begin{align} \label{eq-8.12-f}
\frac{|\nabla v|^{s-2}}{v^q} & \left[|\nabla v|^2 \Delta_f v
-v|\nabla\nabla v|^2 - \frac{v\langle \nabla f , \nabla v \rangle^2}{m-n}\right] \nonumber\\
&\le \frac{|\nabla v|^{s-2}}{v^{q-1}} \left[ \frac{|\nabla v|^2 \Delta_f v}{v} - \frac{(\Delta_f v)^2}{m}\right] \\
&\le \frac{|\nabla v|^{s-2}}{v^{q-1}} 
\left[\frac{m}{4} \frac{|\nabla v|^4}{v^2} -\left( \frac{\sqrt m}{2} \frac{|\nabla v|^2}{v} 
- \frac{\Delta_f v}{\sqrt m} \right)^2\right] \le \frac{m}{4} \frac{|\nabla v|^{s+2}}{v^{q+1}}. \nonumber 
\end{align}
Therefore substituting the above inequality in \eqref{R pq phi equality-f} results in
\begin{align}\label{EQ-eq-10.20-f}
{\mathscr L} ({\mathsf H}^{s,q}_\zeta[v]) 
\le&~ \zeta'(t) \frac{|\nabla v|^s}{v^q} - s \zeta(t) \frac{|\nabla v|^{s-2}}{v^q} 
\left[\frac{1}{2}\partial_t g+(1-p)v{\mathscr Ric}_f^m(g)\right] ( \nabla v, \nabla v) \nonumber\\
&+ \frac{s(1-p)m}{4} \zeta(t) \frac{|\nabla v|^{s+2}}{v^{q+1}}
+ s[q(1-p)-1] \zeta(t) \langle \nabla |\nabla v|^2 , \nabla v\rangle \frac{|\nabla v|^{s-2}}{v^{q}}\nonumber\\
&- q[q(1-p) -p] \zeta(t) \frac{|\nabla v|^{s+2}}{v^{q+1}}
- \frac{s(1-p)}{2v^{q-1}} \zeta(t) \langle \nabla |\nabla v|^{s-2} , \nabla |\nabla v|^2 \rangle\nonumber\\
&- s\zeta (t) \frac{|\nabla v|^{s-2}}{v^q}\langle \nabla v, \nabla \Sigma(t,x,v) \rangle
+ \frac{q \zeta(t)}{v^{q+1}} |\nabla v|^s\Sigma(t,x,v)\nonumber\\
&-\Gamma'(v) \Sigma(t,x,v) - [\Gamma'(v) + (1-p)v\Gamma''(v)]|\nabla v|^2.
\end{align}
Now by using the following identity 
\begin{align}
\left \langle \nabla v, \nabla \left(\frac{|\nabla v|^s}{v^q}\right) \right\rangle 
= \frac{s}{2} \frac{|\nabla v|^{s-2}}{v^q} \langle \nabla v, \nabla |\nabla v|^2 \rangle
- q \frac{|\nabla v|^{s+2}}{v^{q+1}}, 
\end{align} 
and substituting back in \eqref{EQ-eq-10.20-f} we can write 
\begin{align}
{\mathscr L} ({\mathsf H}^{s,q}_\zeta[v]) 
\le&~ \zeta'(t) \frac{|\nabla v|^s}{v^q} - s \zeta(t) \frac{|\nabla v|^{s-2}}{v^q} 
\left[\frac{1}{2}\partial_t g+(1-p) v {\mathscr Ric}_f^m(g)\right] ( \nabla v, \nabla v) \nonumber \\
&+\zeta (t) \frac{s(1-p)m}{4} \frac{|\nabla v|^{s+2}}{v^{q+1}}
+ \zeta(t)[(1-p)q^2 - (2-p)q] \frac{|\nabla v|^{s+2}}{v^{q+1}}\nonumber\\
&+ 2 \zeta (t) [q(1-p)-1] \langle \nabla v, \nabla (|\nabla v|^s/v^q) \rangle 
+ \frac{q \zeta(t)}{v^{q+1}} |\nabla v|^s\Sigma(t,x,v)\nonumber\\
&- \zeta (t)\frac{s(1-p)}{2v^{q-1}} \langle \nabla |\nabla v|^{s-2} , \nabla |\nabla v|^2 \rangle
- s\zeta (t) \frac{|\nabla v|^{s-2}}{v^q}\langle \nabla v, \nabla \Sigma(t,x,v) \rangle\nonumber\\
& - \Gamma'(v) \Sigma(t,x,v) - [\Gamma'(v) + (1-p)v\Gamma''(v)] |\nabla v|^2. 
\end{align}
Next rearranging terms and using the super flow inequality \eqref{eq11c-f} leads to
\begin{align}\label{EQ-2.23}
{\mathscr L} ({\mathsf H}^{s,q}_\zeta[v]) 
\le&~ \zeta'(t) \frac{|\nabla v|^s}{v^q} + s \mathsf k \zeta(t)  \frac{|\nabla v|^{s}}{v^q} 
+\zeta (t)(1-p) \left[q^2 - \frac{2-p}{1-p} q+\frac{sm}{4} \right] \frac{|\nabla v|^{s+2}}{v^{q+1}} \nonumber\\
&+ q\zeta (t) \frac{|\nabla v|^s}{v^q}\frac{\Sigma(t,x,v)}{v} 
+ 2 \zeta (t) [q(1-p)-1] \langle \nabla v, \nabla (|\nabla v|^s/v^q) \rangle \nonumber\\
&- s\zeta (t) \frac{|\nabla v|^{s-2}}{v^q}\langle \nabla v, \nabla \Sigma(t,x,v) \rangle 
- \zeta (t)\frac{s(1-p)}{2v^{q-1}} \langle \nabla |\nabla v|^{s-2} , \nabla |\nabla v|^2 \rangle \nonumber\\
& - \Gamma'(v) \Sigma(t,x,v) - [\Gamma'(v) + (1-p)v\Gamma''(v)] |\nabla v|^2.
\end{align}

Now referring to \eqref{Xpq def-f}, by adding and subtracting suitable terms to the right-hand side of \eqref{EQ-2.23}, we can rewrite this as 
\begin{align}
{\mathscr L} ({\mathsf H}^{s,q}_\zeta[v]) 
\le&~ \zeta'(t) \frac{|\nabla v|^s}{v^q} + s \mathsf k \zeta(t) \frac{|\nabla v|^{s}}{v^q} 
+ 2 [q(1-p)-1] \langle \nabla v, \nabla \overbrace{[\zeta(t)|\nabla v|^s/v^q + \Gamma(v)]}^{{\mathsf H}^{s,q}_\zeta[v]} \rangle \nonumber\\
&+ q\zeta (t)\frac{|\nabla v|^s}{v^{q}}\frac{\Sigma(t,x,v)}{v} 
+\zeta (t)(1-p) \left[q^2 - \frac{2-p}{1-p} q+\frac{sm}{4} \right] \frac{|\nabla v|^{s+2}}{v^{q+1}} \nonumber \\
& - s\zeta(t)\frac{|\nabla v|^{s-2}}{v^q}\langle \nabla v, \nabla \Sigma(t,x,v) \rangle 
- \zeta (t)\frac{s(1-p)}{2v^{q-1}} \langle \nabla |\nabla v|^{s-2} , \nabla |\nabla v|^2 \rangle \nonumber \\
&- 2 [q(1-p)-1]\Gamma'(v) |\nabla v|^2 - \Gamma'(v) \Sigma(t,x,v) \nonumber \\
&- [\Gamma'(v)+(1-p)v\Gamma''(v)]  |\nabla v|^2.  
\end{align}

Finally rearranging terms, using \eqref{Xpq def-f}, and $\langle \nabla v, \nabla \Sigma \rangle 
= \langle \nabla v, \Sigma_x \rangle +  \Sigma_v |\nabla v|^2$ 
along with $\langle \nabla |\nabla v|^{s-2}, \nabla |\nabla v|^2\rangle \ge 0$ when $s \ge 2$ 
gives the desired conclusion. Note that the last inequality follows by writing $V=|\nabla v|^2$, setting 
$\alpha=(s-2)/2$ (recall $s\ge2$) and observing that 
$\langle \nabla V^\alpha, \nabla V \rangle = \alpha V^{\alpha-1} |\nabla V|^2 \ge 0$.
\end{proof}

\begin{remark} \label{q1-q2-remark-fde1} 
The quadratic polynomial $q^2-(2-p)/(1-p)q+sm/4$ appearing in \eqref{EQ-eq-4.4-2.n-2-fde1} 
has non-negative discriminant when $1-1/[\sqrt{sm}-1]\le p<1$. If $q$ is chosen between the 
roots $0<q_1\le q_2$ then this polynomial will be non-positive at $q$. Note also that in the 
closed case in corollaries below the left end-point of the $p$ interval {\it can} be admitted. 
In the general case, due to localisation, we need the strict inequality. 
\end{remark}

\begin{corollary} \label{first-cor-last-fde1}
Let $(M,g,e^{-f} dv_g)$ be a smooth metric measure space with $M$ closed. 
Let $u$ be a positive smooth solution to \eqref{eq11} where $0<1-1/[\sqrt{sm}-1]\le p<1$, $s \ge 2$ and set $v =p/(1-p) u^{p-1}$. 
Assume the metric and potential satisfy the super flow \eqref{eq11c-f} with $\mathsf{k} \ge 0$ and for some $q$ between $q_1$ 
and $q_2$ and $a \in {\mathbb R}$ the following hold: 
\begin{itemize}
\item $\Gamma'(v) \Sigma(t,x,v) \le 0$,
\item $\langle \nabla v, \Sigma_x(t,x,v) \rangle \ge 0$,  
\item $(q/s) \Sigma(t,x,v)/v - \Sigma_v(t,x,v) \le a$, 
\item $[2q(1-p)-1]\Gamma'(v)+(1-p)v\Gamma''(v) \ge 0$.
\end{itemize}
Then for all $x \in M$ and $0<t \le T$ we have 
\begin{equation}
[|\nabla v|^s/v^q] (x,t) 
\le e^{s(\mathsf{k}+a)t} \left\{ \max_{M} 
\left[|\nabla v|^s/v^q + \Gamma(v) \right]_{t=0} - \Gamma(v(x,t)) \right\}. 
\end{equation}
\end{corollary}

\begin{proof}
Using \eqref{EQ-eq-4.4-2.n-2-fde1} and the assumptions on $\Gamma$ and $\Sigma$ in the statement of the corollary, we can write for 
${\mathsf H} = {\mathsf H}^{s,q}_\zeta[v] = \zeta(t) |\nabla v|^s/v^q + \Gamma(v)$
\begin{align}  \label{Eq-7.16-fde1}
{\mathscr L} ({\mathsf H}^{s,q}_\zeta[v]) &- 2 [q(1-p)-1] \langle \nabla v, \nabla {\mathsf H}^{s,q}_\zeta[v]\rangle\\
\le& \left\{ \zeta'(t) + s {\mathsf k} \zeta(t) + \zeta(t) \left[q\frac{\Sigma(t,x,v)}{v} 
- s \Sigma_v(t,x,v) \right] \right\} \frac{|\nabla v|^s}{v^q} \nonumber\\
&+\zeta (t) \left\{ (1-p) \left[q^2 - \frac{2-p}{1-p} q+\frac{sm}{4} \right] \frac{|\nabla v|^{s+2}}{v^{q+1}} 
- s \frac{|\nabla v|^{s-2}}{v^q}\langle \nabla v, \Sigma_x (t,x,v) \rangle \right\} \nonumber\\
&+ \Gamma'(v) \Sigma(t,x,v) - (1-p) \left[\frac{2q(1-p)-1}{1-p} \Gamma'(v)+v\Gamma''(v)\right] |\nabla v|^2 \nonumber \\
\le&~[\zeta'(t) + s ({\mathsf k} + a) \zeta(t)] \frac{|\nabla v|^s}{v^q} 
+\zeta (t)(1-p) \left[q^2 - \frac{2-p}{1-p} q+\frac{sm}{4} \right] \frac{|\nabla v|^{s+2}}{v^{q+1}}. \nonumber 
\end{align}
The function $\zeta(t) = e^{-s(\mathsf{k}+a)t}$ is non-negative, smooth and satisfies $\zeta'+s(\mathsf{k}+a)\zeta=0$. 
Thus substituting in \eqref{Eq-7.16-fde1} and noting the range of $p$ we have 
\begin{equation}
{\mathscr L} ({\mathsf H}^{s,q}_\zeta[v]) - 2 [q(1-p)-1] \langle \nabla v, \nabla {\mathsf H}^{s,q}_\zeta[v]\rangle \le 0.
\end{equation} 
The assertion is now a consequence of the weak maximum principle giving 
\begin{align}
{\mathsf H}^{s,q}_\zeta[v](x,t) &= e^{-s(\mathsf{k}+a)t} [|\nabla v|^s/v^q](x,t) + \Gamma(v(x,t)) \nonumber \\
&\le \max_M [|\nabla v|^s/v^q + \Gamma(v)]_{t=0} = \max_M {\mathsf H}^{s,q}_\zeta[v](x,0).
\end{align}
The proof is thus complete. 
\end{proof}

\begin{corollary} \label{second-cor-last-fde1}
Let $(M,g,e^{-f} dv_g)$ be a smooth metric measure space with $M$ closed. 
Let $u$ be a positive smooth solution to \eqref{eq11}, $0<p_c \le p<1$ and set $v =p/(1-p) u^{p-1}$. 
Assume the metric and potential satisfy the super flow \eqref{eq11c-f} with $\mathsf{k} \ge 0$ and 
for some $1<q<1/(1-p)$ the following hold: 
\begin{itemize}
\item $\Sigma(t,x,v) \ge 0$,
\item $\langle \nabla v, \Sigma_x(t,x,v) \rangle \ge 0$, 
\item $2\Sigma_v(t,x,v) - q\Sigma(t,x,v)/v \ge 0$. 
\end{itemize}
Then for $x \in M$ and $0<t \le T$ we have 
\begin{equation}
t \frac{|\nabla v|^2}{v^q}(x,t) \le \frac{1+2{\mathsf k}t}{(1-q)[(1-p)q-1]} \left[ \max_M v^{1-q}(x,0) - v^{1-q}(x,t) \right]. 
\end{equation}
\end{corollary}

\begin{remark} \label{remark-after-second-cor-last-fde1}
Corollary \ref {second-cor-last-fde1} can be extended to the range $0<1-1/[\sqrt{2m}-1] \le p<p_c$ by 
requiring $q$ to lie in $[q_1, q_2]=[q_1, q_2] \cap (1, 1/(1-p))$. For $p_c\le p<1$ it is easily seen 
that $(1,1/(1-p)) \subset (q_1, q_2)$ and so $(1,1/(1-p))=[q_1, q_2] \cap (1, 1/(1-p))$. In contrast in the sub-critical range 
$0<1-1/[\sqrt{2m}-1] \le p<p_c$ we have $[q_1, q_2] \subset (1,1/(1-p))$. Note that $p \ge p_c=1-2/m$ is equivalent to 
both conditions $q_1 \le 1$ and $1/(1-p) \le q_2$. See also Remark \ref{q1-q2-remark-fde1} regarding the reason for 
$p \ge 1-1/[\sqrt{2m}-1]$. 
\end{remark}

\begin{proof}
For ${\mathsf H}[v] = \zeta(t) |\nabla v|^2/v^q + cv^{1-q}$ where $s=2$, $\Gamma(v)=cv^{1-q}$ and $c>0$ is to be determined 
below, we can write using \eqref{EQ-eq-4.4-2.n-2-fde1} in Lemma \ref{ERAs-PandS-group},
\begin{align} \label{Eq-7.16-fde1-ddd}
{\mathscr L} ({\mathsf H}) - 2 &[q(1-p)-1] \langle \nabla v, \nabla {\mathsf H}\rangle \nonumber \\
\le& \left\{ \zeta'(t) + 2 {\mathsf k} \zeta(t) + \zeta(t) \left[q\frac{\Sigma(t,x,v)}{v} - 2 \Sigma_v(t,x,v) \right] \right\} \frac{|\nabla v|^2}{v^q} \nonumber\\
&+\zeta (t)(1-p) \left[q^2 - \frac{2-p}{1-p} q+\frac{m}{2} \right] \frac{|\nabla v|^4}{v^{q+1}} 
- 2\frac{\zeta (t)}{v^q}\langle \nabla v, \Sigma_x (t,x,v) \rangle \nonumber\\
&- c(1-q) [q(1-p)-1] \frac{|\nabla v|^2}{v^q} + c(1-q) \frac{\Sigma(t,x,v)}{v^q}.
\end{align}
Here we have made use of $\Gamma'(v) \Sigma(t,x,v)=c(1-q) \Sigma(t,x,v)/v^q$ and 
\begin{equation}
(1-p) \left[\frac{2q(1-p)-1}{1-p} \Gamma'(v)+v\Gamma''(v)\right] =  c(1-q) [q(1-p)-1] \frac{1}{v^q}.
\end{equation}
Setting $c=\{(1-q)[q(1-p)-1]\}^{-1}>0$ [note $1<q<1/(1-p)$], the assumptions give 
\begin{align} 
{\rm RHS \, } \eqref {Eq-7.16-fde1-ddd}
\le [\zeta'(t) + 2 {\mathsf k} \zeta(t) -1] \frac{|\nabla v|^2}{v^q} 
+\zeta (t)(1-p) \left[q^2 - \frac{2-p}{1-p} q+\frac{m}{2} \right] \frac{|\nabla v|^4}{v^{q+1}}. \nonumber
\end{align}

Now when $\mathsf{k} \ge 0$ by taking $\zeta(t) = t/(1+2\mathsf{k} t)$ we have 
$(\zeta' + 2 \mathsf{k} \zeta -1) \le 0$. Thus since $p \ge 1-1/[\sqrt{2m}-1]$ 
({\it see} Remark \ref{remark-after-second-cor-last-fde1}) we obtain 
${\mathscr L}({\mathsf H}) - 2 [q(1-p)-1] \langle \nabla v, \nabla {\mathsf H} \rangle \le 0$. 
The conclusion now follows by an application of the weak maximum principle.  
\end{proof}

\section {A Hamilton-Souplet-Zhang estimate in the range $p_0(m)< p < 1$}
\label{sec7}

In this section we present another set of gradient estimates for positive solutions to \eqref{eq11}. 
The range is $p_0< p<1$ where $p_0=p_0(m)>0$ is the larger of the values $1/2$ 
and $1-1/\sqrt{m-1}$ ($m \ge 2)$, that is, $p_0=1/2$ when $2 \le m \le 5$ and 
$p_0=1-1/\sqrt{m-1}$ when $m \ge 5$. Thus we always have $p_0(m) \ge 1/2$, and when $m \ge 4$ 
we have $p_0 \le p_c=1-2/m$. Hence for $m \ge 4$, $(p_0,1)$ is larger than $(p_c,1)$ and 
extends into the sub-critical range. Throughout we use the notation 
\begin{align} \label{Sigma-Star-N-definition}
\Sigma^\star(t, x, v) = (p -1/2) v^{1-\frac{1}{p-1/2}} \mathscr N (t,x,v^{\frac{1}{p-1/2}}).
\end{align}

\begin{theorem} \label{thm-6.1-EQ-FDE}
Let $(M, g, d\mu)$ be a complete smooth metric measure space with $d\mu=e^{-f} dv_g$. 
Suppose that the metric and potential are time dependent, of class $\mathscr{C}^2$ and that for suitable constants 
$k, h \ge 0$ and $m \ge n$ satisfy ${\mathscr Ric}_f^m (g) \ge -(m-1)kg$, $\partial_t g \ge -2h g$ in 
the compact space-time cylinder $Q_{R,T}$ with $R, T >0$. Let $u$ be a positive solution 
to \eqref{eq11} with $p_0< p < 1$ and let $v = u^{p-1/2}$ and $M= \sup_{Q_{R,T}} v$. 
Then there exists $C=C(p,m)>0$ such that for every $(x,t)$ in $Q_{R/2,T}$ with $t>t_0-T$ we have 
\begin{align} \label{eq-7.2-FDE-1/2} 
|\nabla v|(x,t)
\le C \left \{ \begin {array}{ll}
\sqrt h M^\frac{p}{2p-1} + \left[ \dfrac{k^{1/4}}{\sqrt R} + \dfrac{1}{R} + \sqrt k \right] M
\\
\\
+\sup_{Q_{R, T}} \left\{v^{\frac{2p}{3(2p-1)}} |\Sigma^\star_x(t,x,v)|^{1/3}\right\} 
\\
\\
+ \dfrac{M^\frac{p}{2p-1}}{\sqrt {t-t_0+T}} 
+ \sup_{Q_{R, T}}\left\{v^{\frac{p}{2p-1}} [\Sigma^\star_v(t,x,v)]_+^{1/2}\right\}
\end{array}
\right\}.
\end{align}
\end{theorem}

\begin{remark} \label{Sigma-Star-N-u-v-remark}
Since $u$ and $v$ are related to one-another via $u=v^{1/(p-1/2)}$ we can write 
\eqref{Sigma-Star-N-definition} as $\Sigma^\star(t,x,v) = (p-1/2) u^{(2p-3)/2} {\mathscr N}(t,x,u)$. 
\end{remark}

Subject to global bounds we have the following global estimate (in space) that results from passing to the limit 
$R \to \infty$ in \eqref{eq-7.2-FDE-1/2}.

\begin{theorem} \label{thm-6.1-EQ-FDE-global}
Let $(M, g, d\mu)$ be a complete smooth metric measure space with $d\mu=e^{-f} dv_g$. 
Suppose that the metric and potential are time dependent, of class $\mathscr{C}^2$ and 
that for  suitable constants $k, h \ge 0$ and $m \ge n$ satisfy 
${\mathscr Ric}_f^m (g) \ge -(m-1)kg$, 
$\partial_t g \ge -2 h g$ on $M \times [t_0-T, t_0]$. Let $u$ be a positive solution 
to \eqref{eq11} with $p_0< p < 1$ and let $v = u^{p-1/2}$ and $M=\sup v$. 
Then there exists $C=C(p,m)>0$ such that for every $x \in M$ and $t_0-T<t <t_0$ we have 
\begin{align}
|\nabla v|(x,t)
\le C \left \{ \begin {array}{ll}
\left[ \sqrt h + \dfrac{1}{\sqrt {t-t_0+T}} \right] M^\frac{p}{2p-1} + \sqrt k M 
\\
\\
+ \sup_{M \times [t_0-T, t_0]}\left\{v^{\frac{p}{2p-1}}  [\Sigma^\star_v(t,x,v)]_+^{1/2}\right\}
\\
\\
+\sup_{M \times [t_0-T, t_0]} \left\{v^{\frac{2p}{3(2p-1)}} |\Sigma^\star_x(t,x,v)|^{1/3}\right\} 
\end{array}
\right\} 
\end{align}
\end{theorem}

In the static case $\partial_t g \equiv 0$ and $\partial_t f \equiv 0$ by taking $h=0$ we obtain 
the following local version of the estimate. The global version follows by passing to the limit $R\to\infty$.

\begin{theorem} \label{thm-6.1-EQ-FDE-static}
Let $(M, g, d\mu)$ be a complete smooth metric measure space with $d\mu=e^{-f} dv_g$ and  
${\mathscr Ric}_f^m (g) \ge -(m-1)kg$ in ${\mathscr B}_R$ for some $k \ge 0$, $m \ge n$ and $R>0$. 
Let $u$ be a positive solution to \eqref{eq11} with $p_0< p < 1$, $v = u^{p-1/2}$ and 
$M= \sup_{Q_{R,T}} v$. Then there exists $C=C(p,m)>0$ such that for $(x,t)$ in 
$Q_{R/2,T}$ with $t>t_0-T$ we have 
\begin{align} \label{thm-6.1-EQ-FDE-static-equation}
|\nabla v|(x,t)
\le C \left \{ \begin {array}{ll}
\dfrac{M^\frac{p}{2p-1}}{\sqrt {t-t_0+T}} 
+\sup_{Q_{R, T}} \left\{v^{\frac{2p}{3(2p-1)}} |\Sigma^\star_x(t,x,v)|^{1/3}\right\} 
\\
\\
+ \left[ \dfrac{k^{1/4}}{\sqrt R} + \dfrac{1}{R} + \sqrt k \right] M
+ \sup_{Q_{R, T}}\left\{v^{\frac{p}{2p-1}} [\Sigma^\star_v(t,x,v)]_+^{1/2}\right\}
\end{array}
\right\}.
\end{align}
\end{theorem}

We end this section with an application of the estimates above to parabolic Liouville-type theorems. 
Compare with Theorem \ref{ancient-1}.

\begin{theorem} \label{ancient-2}
Let $(M, g, d\mu)$ be a complete smooth metric measure space with $d\mu=e^{-f}dv_g$ and 
${\mathscr Ric}^m_f(g) \ge 0$. Assume $(3-2p) {\mathscr N}(u)/u - 2 {\mathscr N}_u(u) \ge 0$ 
for all $u>0$ where $p_0<p<1$. Then any positive ancient solution to the nonlinear 
fast diffusion equation 
\begin{equation} \label{ancient-equation-2}
\square_p u = \partial_t u - \Delta_f u^p = {\mathscr N}(u(x,t)), 
\end{equation}
satisfying the growth at infinity $u(x,t) = o([\varrho(x) + \sqrt{|t|}]^{2/p})$ 
must be spatially constant. If, in addition, ${\mathscr N}(u) \ge a$ for some $a>0$ and 
all $u>0$ then \eqref{ancient-equation-2} admits no such ancient solutions. 
\end{theorem}

\section{Evolution inequalities ${\bf II}$: ${\mathscr L} = \partial_t- pv^{2(p-1)/(2p-1)} \Delta_f$ 
acting on $w = |\nabla v|^2$} 
\label{sec8}

In this section we derive evolution inequalities for $w = |\nabla v|^2$ where $v$ relates 
to $u$ through $v=u^{p-1/2}$ and ${\mathscr L} = \partial_t- pv^{2(p-1)/(2p-1)} \Delta_f$.

\begin{lemma}\label{FDE-Lem.6.1}
Let $u$ be a positive solution to \eqref{eq11} with $1/2<p<1$ and let $v = u^{p-1/2}$ and 
$\Sigma^\star=\Sigma^\star(t,x,v)$ be as in \eqref{Sigma-Star-N-definition}. Then $v$ satisfies the evolution equation 
\begin{align}\label{FDE-6.1}
{\mathscr L}[v] = \left[ \partial_t- pv^{2(p-1)/(2p-1)} \Delta_f \right] v 
= \frac{p}{2p-1} v^{\frac{1}{1-2p}} |\nabla v|^2 + \Sigma^\star(t,x,v). 
\end{align}
\end{lemma}

\begin {proof}
Since here $u = v^{\frac{1}{p-1/2}}$, a  straightforward differentiation leads to the identities
\begin{align}
\partial_t u =\partial_t v^{\frac{1}{p-1/2}} = \frac{2v^{\frac{1}{p-1/2}-1}}{2p-1} \partial_t v, \qquad 
\Delta_f u^p 
=  \frac{2pv^{\frac{p}{p-1/2}-1}}{2p-1} \left[\Delta_f v + \frac{1}{2p-1} \frac{|\nabla v|^2}{v}\right].
\end{align}
Using these ingredients in \eqref{eq11} then results in  
\begin{align}
\square_p u = \partial_t u - \Delta_f u^p &= \frac{2v^{\frac{1}{p-1/2}-1}}{2p-1} \partial_t v 
-\frac{2pv^{\frac{p}{p-1/2}-1}}{2p-1} \left[\Delta_f v + \frac{|\nabla v |^2}{(2p-1)v} \right] 
= \mathscr N\left(t,x,v^{\frac{1}{p-1/2}}\right) \nonumber
\end{align}
which upon recalling \eqref{Sigma-Star-N-definition} and rearranging terms leads to the desired conclusion. 
\end{proof}

\begin{lemma}\label{FDE-lem-6.2}
Let $(M, g, d\mu)$ be a complete smooth metric measure space with $d\mu=e^{-f} dv_g$. 
Suppose that the metric and potential are time dependent and of class $\mathscr{C}^2$. 
Let $u$ be a positive solution to \eqref{eq11} with $1/2<p<1$ and set $w = |\nabla v|^2$ where $v = u^{p-1/2}$. 
Then $w$ satisfies the evolution inequality 
\begin{align}\label{FDE-eq-6.9}
{\mathscr L}[w] 
\le &~ \frac{2p[(m-1)(p-1)^2-1]}{(2p-1)^2} v^{\frac{2p}{1-2p}} w^2
- [\partial_t g +2p v^{\frac{p-1}{p-1/2}} {\mathscr Ric}_f^m(g)] (\nabla v, \nabla v) \nonumber\\
&+\frac{2p^2v^{\frac{1}{1-2p}}}{2p-1} \langle \nabla v , \nabla w \rangle 
+ 2 \langle \nabla v, \Sigma^\star_x(t,x,v)\rangle +2 \Sigma^\star_v(t,x,v) w. 
\end{align} 
\end{lemma}

\begin{proof}
Firstly, since $w = |\nabla v|^2$ and $\partial_t w = \partial_t (|\nabla v|^2) 
= -[\partial_t g] ( \nabla v, \nabla v) + 2 \langle\nabla v, \nabla \partial_t v \rangle$, a direct calculation and using \eqref {FDE-6.1} gives
\begin{align}\label{FDE-6.8}
\partial_t w
=&-[\partial_t g] ( \nabla v, \nabla v) + 2\left\langle \nabla v, \nabla \left[pv^{\frac{p-1}{p-1/2}} \Delta_f v 
+ \frac{pv^{\frac{1}{1-2p}}}{2p-1} |\nabla v|^2 + \Sigma^\star(t, x, v) \right] \right\rangle  \nonumber\\
=&-[\partial_t g]( \nabla v, \nabla v) + \frac{4p(p -1)}{2p-1} v^{\frac{1}{1-2p}}|\nabla v|^2 \Delta_f v
+2pv^{\frac{p-1}{p-1/2}} \langle \nabla v, \nabla \Delta_f v \rangle \nonumber \\
&+\frac{2pv^{\frac{1}{1-2p}}}{2p-1} \langle \nabla v , \nabla |\nabla v|^2 \rangle 
-\frac{2pv^{\frac{2p}{1-2p}}}{(2p-1)^2} |\nabla v|^4 + 2 \langle \nabla v, \nabla \Sigma^\star(t,x,v)\rangle.
\end{align}
Next, the weighted Bochner-Weitzenb\"ock formula as applied to $v$ describes $\Delta_f w$ as 
\begin{align}\label{EQ-eq-6.9}
\Delta_f w &=\Delta_f |\nabla v|^2 
=2|\nabla \nabla v|^2 +2\langle \nabla v , \nabla \Delta_f v \rangle+ 2{\mathscr Ric}_f^m (\nabla v, \nabla v) 
+ 2 \frac{\langle \nabla f , \nabla v \rangle^2}{m-n}.
\end{align} 
Thus by putting \eqref{FDE-6.8}-\eqref{EQ-eq-6.9} together, rearranging terms and making use of the identities 
$\Delta_f v = \Delta v -\langle \nabla f, \nabla v \rangle$ and $\langle \nabla v, \nabla \Sigma^\star \rangle
=\langle \nabla v, \Sigma^\star_x \rangle + \Sigma^\star_v |\nabla v|^2$ we have 
\begin{align}\label{FDE-6.11}
[\partial_t-pv^{\frac{p-1}{p-1/2}} \Delta_f]w 
=& \overbrace{-2p v^{\frac{p-1}{p-1/2}} |\nabla \nabla v|^2 
-\frac{4p(p-1)}{2p-1} v^{\frac{1}{1-2p}} [\nabla \nabla v(\nabla v, \nabla v) 
- |\nabla v|^2 \Delta v]}^{=-2p v^{1/(1-2p)}{\bf I}} \nonumber \\
& +\frac{2p^2v^{\frac{1}{1-2p}}}{2p-1} \langle \nabla v , \nabla |\nabla v|^2 \rangle 
- [\partial_t g + 2p v^{\frac{p-1}{p-1/2}} {\mathscr Ric}_f^m(g)] (\nabla v, \nabla v) \nonumber\\
&-\frac{2pv^{\frac{2p}{1-2p}}}{(2p-1)^2} |\nabla v|^4 
+ 2 \langle \nabla v, \Sigma^\star_x(t,x,v)\rangle +2 \Sigma^\star_v(t,x,v) |\nabla v|^2 \nonumber \\
&\underbrace{-2p v^{\frac{p-1}{p-1/2}} \left[ \frac{\langle \nabla f, \nabla v \rangle^2}{m-n} 
+ \frac{2(p -1)}{2p-1} |\nabla v|^2 \frac{\langle \nabla f, \nabla v \rangle}{v} \right]}_{=-2p v^{(p-1)/(p-1/2)}{\bf II}}. 
\end{align}

To proceed further with \eqref{FDE-6.11} we need to make use of two inequalities that we now describe. 
As for the first, denoting by ${\bf I}$ the sum of the expressions on the first line in \eqref {FDE-6.11} 
modulo the factor $-2p v^{1/(1-2p)}$, and by $a=(p-1)/(2p-1)$, we have 
\begin{align} \label{justification-hessian-v-inequality}
{\bf I} &= v |\nabla \nabla v|^2 + 2a [\nabla \nabla v(\nabla v, \nabla v) - |\nabla v|^2 \Delta v] \nonumber \\
&= v \left[ |\nabla \nabla v| + a \frac{\nabla \nabla v (\nabla v, \nabla v) - |\nabla v|^2 \Delta v}{v|\nabla \nabla v|} \right]^2 
- a^2 v \left[ \frac{\nabla \nabla v (\nabla v, \nabla v) - |\nabla v|^2 \Delta v}{v|\nabla \nabla v|} \right]^2 \nonumber \\
&\ge - a^2 v\left[ \frac{\nabla \nabla v (\nabla v, \nabla v)}{v|\nabla \nabla v|} 
- \frac{|\nabla v|^2 \Delta v}{v|\nabla \nabla v|} \right]^2 
=  - a^2 \frac{|\nabla v|^4}{v} \left[ \frac{\nabla \nabla v}{|\nabla \nabla v|} 
\left( \frac{\nabla v}{|\nabla v|}, \frac{\nabla v}{|\nabla v|} \right) - \frac{\Delta v}{|\nabla \nabla v|} \right]^2 \nonumber \\
&\ge - (n-1)a^2 |\nabla v|^4/v.
\end{align}
Here in concluding the last inequality we have used the (symmetric) matrix variational identity ({\it see}, e.g., \cite{XXu} Lemma A1, pp.~1419)
\begin{align} \label{matrix-variational-inequality}
\left[ \frac{A(\xi,\xi)}{|A|} - \frac{{\rm tr} A}{|A|} \right]^2 
\le \max_{A \in \mathbf{Sym}_n \atop{|\xi|=1}} \left[ \frac{A(\xi,\xi)}{|A|} - \frac{{\rm tr} A}{|A|} \right]^2 = (n-1).
\end{align}
Thus in summary we have shown, 
\footnote{Note that this inequality is trivially true at points where either $\nabla \nabla v$ or $\nabla v$ vanishes. 
Hence in \eqref {justification-hessian-v-inequality} there is no loss of generality in assuming that neither of these quantities is zero.}
\begin{align}
v |\nabla \nabla v|^2 + \frac{2(p-1)}{2p-1} [\nabla \nabla v(\nabla v, \nabla v) - |\nabla v|^2 \Delta v] + \frac{(n-1)(p-1)^2}{(2p -1)^2} \frac{|\nabla v|^4}{v} \ge 0.
\end{align}
As for the second inequality, denoting by ${\bf II}$ the expression on the last line in \eqref{FDE-6.11} modulo the factor 
$-2p v^{(p-1)/(p-1/2)}$, we have 
\begin{align}
{\bf II} =&~\left[\frac{\langle \nabla f , \nabla v \rangle^2}{m-n} 
+\frac{2(p -1)}{2p-1} |\nabla v|^2 \frac{\langle \nabla f, \nabla v \rangle}{v}\right] 
= \left[ \frac{\langle \nabla f , \nabla v \rangle}{\sqrt{m-n}}
+ \frac{p-1}{2p-1}|\nabla v|^2 \frac{\sqrt{m-n}}{v}\right]^2 \nonumber \\
&- \frac{(m-n)(p-1)^2}{(2p-1)^2} \frac{|\nabla v|^4}{v^2} 
\ge - \frac{(m-n)(p-1)^2}{(2p-1)^2} \frac{|\nabla v|^4}{v^2}.
\end{align}

Returning to \eqref{FDE-6.11} and using the above two inequalities leads after some calculations to 
\begin{align}
[\partial_t - pv^{\frac{p-1}{p-1/2}} \Delta_f]w
\le &~ \frac{2p(m-1)(p-1)^2}{(2p -1)^2} v^{\frac{2p}{1-2p}} |\nabla v|^4 
+\frac{2p^2v^{\frac{1}{1-2p}}}{2p-1} \langle \nabla v , \nabla |\nabla v|^2 \rangle \nonumber\\
&-\frac{2pv^{\frac{2p}{1-2p}}}{(2p-1)^2} |\nabla v|^4
-[\partial_t g+2p v^{\frac{p-1}{p-1/2}} {\mathscr Ric}_f^m(g)] (\nabla v, \nabla v)\nonumber\\
&+ 2 \langle \nabla v, \Sigma^\star_x(t,x,v)\rangle +2 \Sigma^\star_v(t,x,v) |\nabla v|^2. 
\end{align}
The conclusion now follows by making the substitution $w=|\nabla v|^2$ and rearranging.
\end{proof}

\section{Localisation in space-time cylinders $Q_{R,T}$ and proof of Theorem \ref{thm-6.1-EQ-FDE}}
\label{sec9}

Having derived the necessary evolution inequalities, we now proceed to finalising the proof of the local estimate 
in Theorem \ref{thm-6.1-EQ-FDE}. The idea is to start by considering the localised function $\eta w$ where 
$\eta$ is the cylindrical cut-off function in \eqref{cut-off def} and then relate the evolutions of $\eta w$ and $w$ 
under $\mathscr L$ by using an identity similar to \eqref{eta-w-Lpv-equation}. Indeed it is easily seen that here 
for the evolution operator $\mathscr L =\partial_t - pv^{(p-1)/(p-1/2)} \Delta_f$ we have 
\begin{align}
{\mathscr L}[\eta w] = w{\mathscr L}[\eta]
- 2pv^{\frac{p-1}{p-1/2}} [\langle\nabla \eta ,\nabla (\eta w) \rangle- |\nabla \eta|^2 w]/\eta
+\eta {\mathscr L}[w].
\end{align}
Using the bound on ${\mathscr L}[w]$ from \eqref{FDE-eq-6.9} in Lemma \ref{FDE-lem-6.2}
and the relation $\eta \langle \nabla v, \nabla w \rangle = \langle \nabla v, \nabla (\eta w) \rangle - w \langle \nabla v, \nabla \eta \rangle$ 
we can then deduce that 
\begin{align}\label{EQ-eq8.18}
{\mathscr L}[\eta w] \le &~
w {\mathscr L}[\eta]
- 2pv^{\frac{p-1}{p-1/2}} \left\langle \frac{\nabla \eta}{\eta} , \nabla (\eta w) \right\rangle
+ 2pv^{\frac{p-1}{p-1/2}} \frac{|\nabla \eta|^2}{\eta} w \\
&+\frac{2p[(m-1)(p-1)^2-1]}{(2p-1)^2} v^{\frac{2p}{1-2p}} \eta w^2
+\frac{2p^2v^{\frac{1}{1-2p}}}{2p-1} \langle \nabla v , \nabla (\eta w) \rangle\nonumber\\
&-\frac{2p^2v^{\frac{1}{1-2p}}}{2p-1} w \langle \nabla v , \nabla \eta \rangle
+2 [p v^{\frac{p-1}{p-1/2}} (m-1) k+h] \eta w\nonumber\\
&+ 2 \eta \langle \nabla v, \Sigma^\star_x(t,x,v)\rangle +2 \Sigma^\star_v(t,x,v) \eta w.\nonumber
\end{align}

Assume now that the localised function $\eta w$ achieves its maximum over the compact set 
$\{(x,t)|d(x,x_0, t) \le R, t_0-T \le t \le \tau\} \subset M \times [t_0-T,t_0]$ at the point $(x_1,t_1)$. We assume that $(\eta w)(x_1, t_1) >0$ 
as otherwise the conclusion of the theorem is true with $w(x, \tau) \le 0$ for all $d(x, x_0, \tau) \le R/2$. In particular it follows from this 
that $t_1>t_0-T$ and subsequently at the point $(x_1,t_1)$  we have 
\begin{itemize}
\item $\partial_t (\eta w) \ge 0$, 
\item $\nabla(\eta w) =0$, 
\item $\Delta_f(\eta w) \le 0$.
\end{itemize} 
As a result from the above it also follows that at the point $(x_1,t_1)$  we have
\begin{equation}
\mathscr L [\eta w] = [\partial_t - pv^{\frac{p-1}{p-1/2}} \Delta_f] (\eta w) \ge 0.
\end{equation} 
Using this inequality we can rewrite  \eqref{EQ-eq8.18} after a rearrangement of terms as 
\begin{align}
\frac{2p[1-(m-1)(p-1)^2]}{(2p-1)^2} v^{\frac{2p}{1-2p}} \eta w^2 \le
&~ w {\mathscr L}[\eta]
+2 [p v^{\frac{p-1}{p-1/2}} (m-1) k+h] \eta w\nonumber\\
&-\frac{2p^2v^{\frac{1}{1-2p}}}{2p-1} w \langle \nabla v , \nabla \eta \rangle
+ 2pv^{\frac{p-1}{p-1/2}} \frac{|\nabla \eta|^2}{\eta} w\nonumber\\
&+ 2 \eta \langle \nabla v, \Sigma^\star_x(t,x,v)\rangle +2 \Sigma^\star_v(t,x,v) \eta w.
\end{align}

Denoting $p[1-(m-1)(p-1)^2]/(2p-1)^2 = 1/\gamma>0$ (the positivity being a direct consequence of the imposed range of $p$) 
and multiplying through by $\gamma v^{2p/(2p-1)}$ gives 
\begin{align}\label{EQ-eq-6.29}
2\eta w^2 \le
&~ \gamma v^{\frac{2p}{2p-1}} w {\mathscr L}[\eta]
+2 \gamma [p v^{\frac{p-1}{p-1/2}} (m-1) k+h] v^{\frac{2p}{2p-1}} \eta w\nonumber\\
&-\frac{2 \gamma p^2}{2p-1} v w \langle \nabla v , \nabla \eta \rangle
+ 2 \gamma p v^2 \frac{|\nabla \eta|^2}{\eta} w\nonumber\\
& + 2 \gamma \eta v^{\frac{2p}{2p-1}} \langle \nabla v, \Sigma^\star_x(t,x,v)\rangle 
+2 \gamma v^{\frac{2p}{2p-1}}\Sigma^\star_v(t,x,v) \eta w.
\end{align}

Now we proceed onto bounding the expression on the right-hand side of \eqref{EQ-eq-6.29}.
As this is similar to the proof of Theorem \ref{thm2-FDE} we shall remain brief, focusing mainly 
on the differences. Starting from the first term and arguing as in 
Theorem \ref{thm2-FDE} we have 
\begin{align}
\gamma v^{\frac{2p}{2p-1}} w \mathscr L[\eta] 
=&~ \gamma v^{\frac{2p}{2p-1}} w[\partial_t - pv^{\frac{p-1}{p-1/2}} \Delta_f] \eta 
= \gamma v^{\frac{2p}{2p-1}} w \partial_t \eta - \gamma p v^2 w \Delta_f \eta\nonumber\\
\le&~ \frac{1}{6} \eta w^2
+ C \gamma ^2 \left[\frac{1}{(\tau-t_0+T)^2} + h^2\right](\sup_{Q_{R, T}} v)^{\frac{4p}{2p-1}}\nonumber\\
&+ C \gamma^2p^2 (m-1)^2\left[\frac{1+ k R^2}{R^4}\right] (\sup_{Q_{R, T}} v)^{4}.
\end{align}
Next, for the second and third terms, we can write 
\begin{align}
2 \gamma [p v^{\frac{p-1}{p-1/2}} (m-1) k] v^{\frac{2p}{2p-1}} \eta w
&= 2 \gamma p (m-1) k v^2 \eta w 
\le \frac{1}{6} \eta w^2 + C \gamma ^2 p^2 (m-1)^2 k^2(\sup_{Q_{R, T}} v)^4, \nonumber \\
2 \gamma h v^{\frac{2p}{2p-1}} \eta w &\le
2 \gamma h v^{\frac{2p}{2p-1}} \eta^{1/2} w 
\le \frac{1}{6} \eta w^2 + C \gamma ^2 h^2 (\sup_{Q_{R, T}} v)^{\frac{4p}{2p-1}}.
\end{align} 
We can bound the fourth term by writing 
\begin{align}
-\frac{2 \gamma p^2}{2p-1} v w \langle \nabla v , \nabla \eta \rangle 
\le \frac{2 \gamma p^2}{2p-1} vw |\nabla v| |\nabla \eta| 
&\le \frac{2 \gamma p^2}{2p-1} v \eta^{3/4} w^{3/2} \frac{|\nabla \eta|}{\eta^{3/4}}\nonumber\\
\le \frac{1}{6} (\eta ^{3/4} w^{3/2})^{4/3} + \frac{C\gamma^4 p^8}{(2p-1)^4} 
\left[\frac{|\nabla \eta|}{\eta^{3/4}}\right]^4 v^4  
&\le \frac{1}{6} \eta w^2 + \frac{C \gamma^4 p^8}{(2p-1)^4 R^4} (\sup_{Q_{R, T}} v)^4.
\end{align}
In much the same way
\begin{align}
2 \gamma p v^2 \frac{|\nabla \eta|^2}{\eta} w 
\le \frac{1}{6} \eta w^2 + \frac{C \gamma^2 p^2}{R^4}(\sup_{Q_{R, T}} v)^4.
\end{align}

Finally for the $\Sigma^\star(t,x,v)$ terms, upon recalling $0\le \eta \le 1$ and utilising Young's inequality we have
\begin{align}
2 \gamma \eta w v^{\frac{2p}{2p-1}}\Sigma^\star_v(t,x,v) \le \frac{1}{6} \eta w^2 + C \gamma^2 \sup_{Q_{R, T}} 
\left\{v^{\frac{4p}{2p-1}} \left[\Sigma^\star_v(t,x,v) \right]_+^2\right\}.
\end{align}
Furthermore, by using Cauchy-Schwarz and Young inequalities we can write
\begin{align}
2 \gamma \eta v^{\frac{2p}{2p-1}} \langle \nabla v, \Sigma^\star_x(t,x,v) \rangle \le &~
2 \gamma \eta v^{\frac{2p}{2p-1}} |\nabla v| |\Sigma^\star_x(t,x,v)| \nonumber \\
\le &~ \frac{1}{6} \eta w^2+C \gamma^{4/3} \sup_{Q_{R, T}} 
\left\{v^{\frac{8p}{3(2p-1)}} |\Sigma^\star_x(t,x,v)|^{4/3}\right\}. 
\end{align}

Completing the estimate of the individual terms on the right-hand side of \eqref{EQ-eq-6.29}, by inserting these estimates
back and after basic considerations and adjusting the constants at the space-time point $(x_1, t_1)$, 
we arrive at the following bound:
\begin{align}\label{EQ-eq-6.39}
\eta w^2 \le&~C\gamma^2 \left \{ \begin {array}{ll} 
p^2 (m-1)^2\left[\dfrac{1+ k R^2}{R^4}\right]
+\dfrac{p^2}{R^4}
\\
\\
+\dfrac{\gamma^2 p^8}{(2p-1)^4 R^4}
+p^2 (m-1)^2 k^2
\end{array}
\right\}
(\sup_{Q_{R, T}} v)^4 \\
&+C \left \{ \begin {array}{ll} 
\gamma ^{2} \left[\dfrac{1}{(\tau-t_0+T)^2} + 2h^2\right](\sup_{Q_{R, T}} v)^{\frac{4p}{2p-1}}
\\
\\
+\gamma ^{4/3}\sup_{Q_{R, T}} \left\{v^{\frac{8p}{3(2p-1)}} |\Sigma^\star_x(t,x,v)|^{4/3}\right\} 
\\
\\
+ \gamma^{2} \sup_{Q_{R, T}} 
\left\{v^{\frac{4p}{2p-1}} \left[\Sigma^\star_v(t,x,v) \right]_+^2\right\} 
\end{array}
\right\}.\nonumber
\end{align}

Upon adjusting terms, absorbing $\gamma$ into $C$ and recalling 
$M=\sup_{Q_{R,T}} v$ it follows that 
\begin{align}
\eta w^2 \le &~C(\gamma) \left \{ \left( \begin {array}{ll} 
p^2 (m-1)^2\left[\dfrac{1+ k R^2}{R^4}\right]
+\dfrac{p^2}{R^4} 
\\
\\
+\dfrac{p^8}{(2p-1)^4 R^4} + p^2 (m-1)^2 k^2
\end{array}
\right) \right. M^4
\nonumber \\ 
&\qquad \qquad + 
\left. \left( 
\begin{array}{ll}
\left[ \dfrac{1}{(\tau-t_0+T)^2} + 2h^2 \right] M^\frac{4p}{2p-1}
\\
\\
+ \sup_{Q_{R, T}} \left\{v^{\frac{8p}{3(2p-1)}} |\Sigma^\star_x(t,x,v)|^{4/3}\right\} 
\\
\\
+\sup_{Q_{R, T}} \left\{v^{\frac{4p}{2p-1}} [\Sigma^\star_v(t,x,v)]_+^2\right\}
\end{array}
\right) \right\}.
\end{align}

Recalling the maximality of $\eta w$ at $(x_1, t_1)$ along with $\eta \equiv 1$ when $d(x, x_0,t) \le R/2$ and $\tau \le t \le t_0$, 
it follows that $w^2 (x, \tau) = (\eta^2 w^2)(x, \tau) \le (\eta^2 w^2)(x_1,t_1) \leq (\eta w^2)(x_1,t_1)$. 
Hence noting $w = |\nabla v|^2$ where $v = u^{p-1/2}$ and 
absorbing $p$ and $m$ into $C(\gamma)$ gives
\begin{align}
|\nabla v|(x,\tau)
\le C \left \{ \begin {array}{ll}
\left[ \dfrac{k^{1/4}}{\sqrt {R}} +\dfrac{1}{R} +\sqrt k \right] M +\sqrt h M^\frac{p}{2p-1}
\\
\\
+ \dfrac{M^\frac{p}{2p-1}}{\sqrt {\tau-t_0+T}} 
+ \sup_{Q_{R, T}}\left\{v^{\frac{p}{2p-1}} [\Sigma^\star_v(t,x,v)]_+^{1/2}\right\}
\\
\\
+ \sup_{Q_{R, T}} \left\{v^{\frac{2p}{3(2p-1)}} |\Sigma^\star_x(t,x,v)|^{1/3}\right\} 
\end{array}
\right\},
\end{align}
where $C=C(p,m)>0$. Making note of the arbitrariness of $\tau$ in the interval $(t_0-T, t_0]$ gives 
the desired conclusion.
\hfill $\square$

\section{Ancient solutions to $\partial_t u - \Delta_f u^p = {\mathscr N}(u)$ ${\bf II}$: Proof of Theorem \ref{ancient-2}}
\label{sec10}

Let us now present the proof of Theorem \ref{ancient-2}. We first fix a space-time point $(x_0,t_0)$. 
Then with the choices $t=t_0$, $R>0$ and $T=R^2$ it follows from 
\eqref{thm-6.1-EQ-FDE-static-equation} in Theorem \ref{thm-6.1-EQ-FDE-static} (with $k=0$) that 
\begin{align} \label{vx0t0-ancient-2}
|\nabla v|(x_0,t_0)
&\le C \left \{ \begin {array}{ll}
\dfrac{M^\frac{p}{2p-1}}{\sqrt {t-t_0+T}}
+\sup_{Q_{R, T}} \left\{v^{\frac{2p}{3(2p-1)}} |\Sigma^\star_x(t,x,v)|^{1/3}\right\} 
\\
\\
+\left[ \dfrac{k^{1/4}}{\sqrt R} + \dfrac{1}{R} + \sqrt k \right] M 
+ \sup_{Q_{R, T}}\left\{v^{\frac{p}{2p-1}} [\Sigma^\star_v(t,x,v)]_+^{1/2}\right\} 
\end{array}
\right\}
\nonumber \\
&
\le C \left \{ \begin {array}{ll}
\dfrac{M^\frac{p}{2p-1}}{\sqrt T} 
+\sup_{Q_{R, T}} \left\{v^{\frac{2p}{3(2p-1)}} |\Sigma^\star_x(t,x,v)|^{1/3} \right\}
\\
\\
+ \dfrac{M}{R} 
+ \sup_{Q_{R, T}} \left\{v^{\frac{p}{2p-1}}\left[\Sigma^\star_v (t,x,v)\right]_+^{1/2} \right\} 
\end{array}
\right\}.
\end{align}

Now $\Sigma^\star(t,x,v) = (p-1/2) u^{(2p-3)/2} {\mathscr N}(u)$ where 
$u=v^{1/(p-1/2)}$ ({\it see} \eqref {ancient-equation-2} and Remark \ref{Sigma-Star-N-u-v-remark}). 
Hence it is seen that $\Sigma^\star$ does not depend explicitly on $t$ and $x$ and 
so $\Sigma^\star = \Sigma^\star (v)$. In particular $\Sigma^\star_x \equiv 0$ 
and by a direct calculation using 
$\Sigma^\star_v (v) = \Sigma^\star_u (v) \partial_v u$ that
\begin{align} 
\Sigma^\star_v (v) 
&= (p-1/2) [(p-3/2) u^{p-5/2} {\mathscr N} (u)+ u^{p-3/2} {\mathscr N}_u(u)] 
u^{(3-2p)/2} /(p-1/2) \nonumber \\
&= (p-3/2) {\mathscr N}(u)/u + {\mathscr N}_u(u). 
\end{align}
From the assumptions on ${\mathscr N}$ it now follows that $\Sigma^\star_v (v) \le 0$ and as a result  
$[\Sigma^\star_v]_+ \equiv 0$. Next from the growth assumption on $u$ we have  
\begin{equation} \label{M-growth-ancient-2}
v = u^{p-1/2} = \left[ o([\varrho(x) + \sqrt{|t|}]^{2/p}) \right]^{p-1/2} = o([\varrho(x) + \sqrt{|t|}]^{(2p-1)/p}), 
\end{equation} 
and so as a result $M=\sup_{Q_{R,T}} v =o(R^{(2p-1)/p})$. Returning now to 
\eqref{vx0t0-ancient-2}, using the above we have
\begin{align}
|\nabla v| (x_0,t_0) \le C \left[ \frac{M^\frac{p}{2p-1}}{\sqrt T}+\frac{M}{R} \right] 
\le \frac{o(R)}{R} + \frac{o(R^{(2p-1)/p})}{R}.
\end{align}
Passing to the limit $R \nearrow \infty$ and noting that $(2p-1)/p \le 1$ in the range $1/2<p<1$ it follows that $|\nabla v|(x_0,t_0)=0$. The arbitrariness of $(x_0,t_0)$ implies 
$|\nabla v| \equiv 0$ and so $v$ and subsequently $u$ are spatially constant. Hence we have $u=u(t)$. From the 
equation satisfied by $u$ it then follows that $du/dt = {\mathscr N}(u)$. The rest of the proof proceeds as in the  proof 
of Theorem \ref{ancient-1}. \hfill $\square$

\section{The Ricci-Perelman super flow ${\bf II}$: $\partial_t g + 2p v^{2(p-1)/(2p-1)} {\mathscr Ric}^m_f(g) \ge - 2\mathsf{k}g$} 
\label{sec11}

In this final section of the paper we derive further global gradient bounds on positive smooth solutions 
to equation \eqref{eq11}. Here $M$ is closed and the metric and potential are assumed to evolve 
under the super flow [{\it cf.} \eqref{SPR-p-substitute-intro} and \eqref{SP-2-intro}]
\begin{align} \label{eq11c-2.n} 
\frac{1}{2} \dfrac{\partial g}{\partial t} (x,t) + p v^{2(p-1)/(2p-1)}(x,t) {\mathscr Ric}^m_f(g)(x,t) 
\ge - \mathsf{k}g(x,t), \qquad {\mathsf k} \ge 0.
\end{align} 

\begin{lemma}\label{lemma_2.1}
Let $u $ be a positive solution to \eqref{eq11} and $v = u^{p-1/2}$ where $1/2<p<1$. 
For $s \ge 2$ and $\zeta=\zeta(t)$ non-negative and of class ${\mathscr C}^1[0, \infty)$ let
\begin{equation}\label{EQ-eq-6.3.n}
{\mathsf H}^{s}_\zeta[v] = \zeta(t) |\nabla v|^s + \Gamma(v),
\end{equation} 
where $\Gamma=\Gamma(v)$ is of class $\mathscr{C}^2(0, \infty)$. 
Then ${\mathsf H} = {\mathsf H}^{s}_\zeta [v]$ satisfies the evolution identity
\begin{align}\label{identity-2.4}
{\mathscr L} ({\mathsf H}^{s}_\zeta[v]) - 
&\frac{2p^2v^{\frac{1}{1-2p}}}{2p-1} \langle \nabla v , \nabla {\mathsf H}^{s}_\zeta[v] \rangle\nonumber\\
=&~\zeta'(t) |\nabla v|^s - s\zeta (t) |\nabla v|^{s-2} \left[ \frac{1}{2} \partial_t g 
+ pv^{\frac{p-1}{p-1/2}} {\mathscr Ric}_f^m(g) \right] ( \nabla v, \nabla v) \nonumber\\
& - \frac{2 sp(p-1)}{2p-1} \zeta(t) v^{\frac{1}{1-2p}} |\nabla v|^{s-2} \nabla \nabla v (\nabla v, \nabla v)
- \frac{sp \zeta(t)}{(2p-1)^2} v^{\frac{2p}{1-2p}} |\nabla v|^{s+2} \nonumber\\
&- sp\zeta(t) v^{\frac{p-1}{p-1/2}} |\nabla v|^{s-2} |\nabla \nabla v|^2 
+ \frac{2sp(p -1)}{2p-1} \zeta(t) v^{\frac{1}{1-2p}} |\nabla v|^s \Delta v \nonumber\\
&- \frac{s}{2} \zeta(t) pv^{\frac{p-1}{p-1/2}} \langle \nabla |\nabla v|^{s-2} , \nabla |\nabla v|^2 \rangle 
- \frac{2sp(p -1)}{2p-1} \zeta(t) v^{\frac{1}{1-2p}} |\nabla v|^s \langle \nabla f, \nabla v \rangle \nonumber\\
& - sp \zeta(t) v^{\frac{p-1}{p-1/2}} |\nabla v|^{s-2} \langle \nabla f , \nabla v \rangle^2/(m-n)
+ s \zeta(t) |\nabla v|^{s-2}\langle \nabla v, \nabla \Sigma^\star (t,x,v)\rangle \nonumber\\
&+\Gamma'(v) \Sigma^\star (t,x,v) -p v^{\frac{1}{1-2p}} [\Gamma'(v)+v\Gamma''(v)] |\nabla v|^2.
\end{align}
\end{lemma}

\begin{proof}
Referring to \eqref{EQ-eq-6.3.n} it is a straightforward matter to see that  
\begin{align}\label{EQ-eq-6.6.n}
{\mathscr L}({\mathsf H}^{s}_\zeta[v]) = \zeta (t) {\mathscr L}[|\nabla v|^s ] 
+ \zeta'(t) |\nabla v|^s + {\mathscr L}[\Gamma(v)].
\end{align}

We now go forward by evaluating the expressions in the first and third terms on the right-hand side. 
For the first term, ignoring the factor $\zeta(t)$ for the moment, we have  
\begin{align}\label{EQ-eq-6.11.n}
{\mathscr L} [|\nabla v|^s] 
=&~\frac{s}{2} |\nabla v|^{s-2} \left[ \partial_t -pv^{\frac{p-1}{p-1/2}} \Delta_f\right] |\nabla v|^2
- \frac{s}{2}pv^{\frac{p-1}{p-1/2}} \langle \nabla |\nabla v|^{s-2} , \nabla |\nabla v|^2 \rangle,  
\end{align}
where by a direct calculation using Lemma \ref{FDE-Lem.6.1} and \eqref{Bochner-1} it is seen that 
\begin{align}
\left[ \partial_t -pv^{\frac{p-1}{p-1/2}} \Delta_f\right] |\nabla v|^2
= & -[\partial_t g] ( \nabla v, \nabla v) +\frac{4p(p -1)}{2p-1} v^{\frac{1}{1-2p}} |\nabla v|^2 \Delta_f v \\
& +\frac{2pv^{\frac{1}{1-2p}}}{2p-1} \langle \nabla v , \nabla |\nabla v|^2 \rangle 
-\frac{2pv^{\frac{2p}{1-2p}} |\nabla v|^4}{(2p-1)^2} \nonumber\\
&+2 \langle \nabla v, \nabla \Sigma^\star(t,x,v) \rangle-2 pv^{\frac{p-1}{p-1/2}} |\nabla \nabla v|^2\nonumber\\
& -2 pv^{\frac{p-1}{p-1/2}} {\mathscr Ric}_f^m (\nabla v, \nabla v) 
- 2 pv^{\frac{p-1}{p-1/2}} \frac{\langle \nabla f , \nabla v \rangle^2}{m-n}. \nonumber 
\end{align}
Therefore by combining the last two identities it follows that 
\begin{align}\label{EQ-eq-6.13.n}
{\mathscr L}[|\nabla v|^s] 
=&-\frac{s}{2} |\nabla v|^{s-2} [\partial_t g] ( \nabla v, \nabla v) 
+ \frac{2sp(p -1)}{2p-1} v^{\frac{1}{1-2p}} |\nabla v|^s \Delta_f v\nonumber\\
&+ \frac{spv^{\frac{1}{1-2p}}}{2p-1} |\nabla v|^{s-2} \langle \nabla v, \nabla |\nabla v|^2 \rangle
- \frac{spv^{\frac{2p}{1-2p}}}{(2p-1)^2} |\nabla v|^{s+2}\nonumber\\
&+ s |\nabla v|^{s-2}\langle \nabla v, \nabla \Sigma^\star(t,x,v) \rangle
- s pv^{\frac{p-1}{p-1/2}} |\nabla v|^{s-2} |\nabla \nabla v|^2\nonumber\\
& - s pv^{\frac{p-1}{p-1/2}} |\nabla v|^{s-2} {\mathscr Ric}_f^m (\nabla v, \nabla v)
-s pv^{\frac{p-1}{p-1/2}} |\nabla v|^{s-2} \frac{\langle \nabla f , \nabla v \rangle^2}{m-n}\nonumber\\
&-\frac{s}{2} pv^{\frac{p-1}{p-1/2}} \langle \nabla |\nabla v|^{s-2} , \nabla |\nabla v|^2 \rangle.
\end{align}
In a similar way for the third term on the right-hand side of \eqref{EQ-eq-6.6.n} we can write 
\begin{align}\label{EQ-eq-6.14.n}
{\mathscr L}[\Gamma(v)] 
&= \Gamma'(v) \left [\frac{pv^{\frac{1}{1-2p}} |\nabla v|^2}{2p-1} + \Sigma^\star(t, x, v)\right] - pv^{\frac{p-1}{p-1/2}} \Gamma''(v)  |\nabla v|^2, 
\end{align}
where we have used \eqref{FDE-6.1}. Substituting \eqref{EQ-eq-6.13.n} and \eqref{EQ-eq-6.14.n} back into \eqref{EQ-eq-6.6.n} now leads to
\begin{align}\label{EQ-eq-6.15.n}
{\mathscr L} ({\mathsf H}^{s}_\zeta[v]) 
=&~ \zeta'(t) |\nabla v|^s - s\zeta(t) |\nabla v|^{s-2}\left[\frac{1}{2}\partial_t g 
+ pv^{\frac{p-1}{p-1/2}} {\mathscr Ric}_f^m(g)\right] ( \nabla v, \nabla v) \nonumber\\
&+ \frac{2sp(p -1)}{2p-1} \zeta(t) v^{\frac{1}{1-2p}} |\nabla v|^s [\Delta v - \langle \nabla f, \nabla v \rangle]
- sp \zeta(t)v^{\frac{p-1}{p-1/2}} |\nabla v|^{s-2} |\nabla \nabla v|^2\nonumber\\
&+ \frac{sp \zeta(t)}{2p-1} v^{\frac{1}{1-2p}} |\nabla v|^{s-2}\langle \nabla v, \nabla |\nabla v|^2 \rangle 
-\frac{sp \zeta(t)}{(2p-1)^2} v^{\frac{2p}{1-2p}} |\nabla v|^{s+2} \nonumber\\
& - sp\zeta(t)v^{\frac{p-1}{p-1/2}} |\nabla v|^{s-2} \frac{\langle \nabla f, \nabla v \rangle^2}{m-n}
- \frac{sp}{2} \zeta(t) v^{\frac{p-1}{p-1/2}} \langle \nabla |\nabla v|^{s-2} , \nabla |\nabla v|^2 \rangle \nonumber\\
& + s\zeta(t)|\nabla v|^{s-2}\langle \nabla v, \nabla \Sigma^\star(t,x,v) \rangle 
- pv^{\frac{p-1}{p-1/2}} \Gamma''(v)  |\nabla v|^2\nonumber\\
&+ \Gamma'(v) \left [\frac{p v^{\frac{1}{1-2p}} |\nabla v|^2}{2p-1} + \Sigma^\star(t, x, v)\right].
\end{align}
Next, upon rearranging terms and splitting the term $v^{\frac{1}{1-2p}} |\nabla v|^{s-2}\langle \nabla v , \nabla |\nabla v|^2 \rangle$
we have 
\begin{align}
{\mathscr L}({\mathsf H}^{s}_\zeta[v]) 
=&- sp\zeta(t) v^{\frac{p-1}{p-1/2}} |\nabla v|^{s-2} |\nabla \nabla v|^2 
+ \frac{2sp(p -1)}{2p-1} \zeta(t) v^{\frac{1}{1-2p}} |\nabla v|^s \Delta v \nonumber\\
&-\frac{2 sp(p-1)}{2p-1} \zeta(t) v^{\frac{1}{1-2p}} |\nabla v|^{s-2} \nabla \nabla v (\nabla v, \nabla v)
+ \zeta'(t) |\nabla v|^s \nonumber\\
&+\frac{sp^2 \zeta(t)}{2p-1} v^{\frac{1}{1-2p}} |\nabla v|^{s-2} \langle \nabla v , \nabla |\nabla v|^2 \rangle
 -\frac{sp\zeta(t)}{(2p-1)^2} v^{\frac{2p}{1-2p}} |\nabla v|^{s+2} \nonumber\\
&- s\zeta (t) |\nabla v|^{s-2}\left[\frac{1}{2}\partial_t g + pv^{\frac{p-1}{p-1/2}} 
{\mathscr Ric}_f^m(g)\right] (\nabla v, \nabla v) \nonumber\\
&- \frac{sp}{2} \zeta(t) v^{\frac{p-1}{p-1/2}} \langle \nabla |\nabla v|^{s-2} , \nabla |\nabla v|^2 \rangle 
- \frac{2sp(p -1)}{2p-1} \zeta(t) v^{\frac{1}{1-2p}} |\nabla v|^s \langle \nabla f, \nabla v \rangle \nonumber\\
&- sp\zeta(t)v^{\frac{p-1}{p-1/2}} |\nabla v|^{s-2} \frac{\langle \nabla f, \nabla v \rangle^2}{m-n}
+ s\zeta(t)|\nabla v|^{s-2}\langle \nabla v, \nabla \Sigma^\star(t,x,v) \rangle \nonumber\\
&- pv^{\frac{p-1}{p-1/2}} \Gamma''(v)  |\nabla v|^2
+ \Gamma'(v) \left [\frac{pv^{\frac{1}{1-2p}} |\nabla v|^2}{2p-1} + \Sigma^\star(t, x, v)\right].
\end{align}
Taking advantage of the relation
$2\langle \nabla v , \nabla |\nabla v|^s \rangle = s |\nabla v|^{s-2} \langle \nabla v, |\nabla v|^2 \rangle$ then gives
\begin{align}
{\mathscr L}({\mathsf H}^{s}_\zeta[v]) - 
&\frac{2p^2v^{\frac{1}{1-2p}}}{2p-1} \langle \nabla v , \nabla {\mathsf H}^{s}_\zeta[v] \rangle\nonumber\\
=&- sp\zeta(t)v^{\frac{p-1}{p-1/2}} |\nabla v|^{s-2} |\nabla \nabla v|^2 
+ \frac{2sp(p -1)}{2p-1} \zeta(t)v^{\frac{1}{1-2p}} |\nabla v|^s \Delta v \nonumber\\
& -\frac{2 sp(p-1)}{2p-1} \zeta(t)v^{\frac{1}{1-2p}} |\nabla v|^{s-2} \nabla \nabla v (\nabla v, \nabla v)
-\frac{sp \zeta(t)}{(2p-1)^2} v^{\frac{2p}{1-2p}} |\nabla v|^{s+2} \nonumber\\
&+ \zeta'(t) |\nabla v|^s - s\zeta (t)|\nabla v|^{s-2}\left[\frac{1}{2}\partial_t g 
+ pv^{\frac{p-1}{p-1/2}} {\mathscr Ric}_f^m(g)\right] (\nabla v, \nabla v) \nonumber\\
&- \frac{sp}{2} \zeta(t) v^{\frac{p-1}{p-1/2}} \langle \nabla |\nabla v|^{s-2} , \nabla |\nabla v|^2 \rangle 
- s\zeta(t) \frac{2p(p -1)}{2p-1} v^{\frac{1}{1-2p}} |\nabla v|^s \langle \nabla f, \nabla v \rangle \nonumber\\
& - sp\zeta(t)v^{\frac{p-1}{p-1/2}} |\nabla v|^{s-2} \frac{\langle \nabla f, \nabla v \rangle^2}{m-n}
+ s\zeta(t) |\nabla v|^{s-2}\langle \nabla v, \nabla \Sigma^\star(t,x,v) \rangle \nonumber\\
&+ \Gamma'(v) \Sigma^\star(t, x, v) - p v^{\frac{1}{1-2p}} \Gamma'(v) |\nabla v|^2
- pv^{\frac{p-1}{p-1/2}} \Gamma''(v)  |\nabla v|^2,
\end{align}
which at once gives the desired conclusion.
\end{proof}

\begin{lemma} \label{evolution Rpq-2.n-inequality-2}
Under the assumptions of Lemma $\ref{lemma_2.1}$, if the metric and potential evolve under 
the super flow \eqref{eq11c-2.n}, then ${\mathsf H} = {\mathsf H}^{s}_\zeta [v]$ satisfies the evolution inequality 
\begin{align} \label{EQ-eq-4.4-2.n-2}
{\mathscr L} & ({\mathsf H}^{s}_\zeta[v]) - 
\frac{2p^2v^{\frac{1}{1-2p}}}{2p-1} \langle \nabla v , \nabla {\mathsf H}^{s}_\zeta[v] \rangle \\
\le & [\zeta'(t)+ s \mathsf k \zeta(t) + s \zeta(t) \Sigma_v^\star(t,x,v)] |\nabla v|^s
+ \frac{s p[(p-1)^2(m-1)-1]}{(2p-1)^2} \zeta(t) v^{\frac{2p}{1-2p}} |\nabla v|^{s+2} \nonumber\\
&+ s \zeta(t) |\nabla v|^{s-2}\langle \nabla v,  \Sigma^\star_x(t,x,v) \rangle 
+ \Gamma'(v) \Sigma^\star(t,x,v) - pv^{\frac{1}{1-2p}} [\Gamma'(v)+v\Gamma''(v)]  |\nabla v|^2. \nonumber 
\end{align}
\end{lemma}

\begin{proof}
We start with identity \eqref{identity-2.4} in Lemma \ref{lemma_2.1}. 
Using $A=\nabla \nabla v$ with
${\rm tr} A=\Delta v$ and $\xi = \nabla v/ |\nabla v|$ 
we can rewrite this upon substitution and a rearrangement of terms in the form 
\begin{align}
{\mathscr L} &({\mathsf H}^{s}_\zeta[v]) - 
\frac{2p^2v^{\frac{1}{1-2p}}}{2p-1} \langle \nabla v , \nabla {\mathsf H}^{s}_\zeta[v] \rangle\nonumber\\
=&- sp\zeta(t) v^{\frac{p-1}{p-1/2}} |\nabla v|^{s-2} |A|^2 
- \frac{2sp(p -1)}{2p-1} \zeta(t) v^{\frac{1}{1-2p}} \left[ \frac{A (\xi, \xi)}{|A|}- \frac{{\rm tr} A}{|A|}\right]|\nabla v|^s |A| \nonumber\\
& -\frac{sp\zeta(t)}{(2p-1)^2} v^{\frac{2p}{1-2p}} |\nabla v|^{s+2} 
 -s\zeta (t)|\nabla v|^{s-2}\left[\frac{1}{2}\partial_t g + pv^{\frac{p-1}{p-1/2}} {\mathscr Ric}_f^m(g)\right] ( \nabla v, \nabla v)\nonumber\\
&+ \zeta'(t) |\nabla v|^s-\frac{sp}{2} \zeta(t) v^{\frac{p-1}{p-1/2}} \langle \nabla |\nabla v|^{s-2} , \nabla |\nabla v|^2 \rangle 
+ s\zeta(t) |\nabla v|^{s-2}\langle \nabla v, \nabla \Sigma^\star(t,x,v) \rangle\nonumber\\
& - sp\zeta(t) v^{\frac{p-1}{p-1/2}} |\nabla v|^{s-2} \left[\frac{\langle \nabla f, \nabla v \rangle^2}{(m-n)}
+ \frac{2(p -1)}{2p-1} |\nabla v|^2 \frac{\langle \nabla f, \nabla v \rangle}{v}\right] \nonumber\\
&+\Gamma'(v) \Sigma^\star (t,x,v) -p v^{\frac{1}{1-2p}} [\Gamma'(v)+v\Gamma''(v)] |\nabla v|^2.
\end{align}
Upon forming a complete square for the matrix term it then follows that
\begin{align}
{\mathscr L} &({\mathsf H}^{s}_\zeta[v]) - 
\frac{2p^2 v^{\frac{1}{1-2p}}}{2p-1} \langle \nabla v , \nabla {\mathsf H}^{s}_\zeta[v] \rangle\nonumber\\
=&- sp\zeta(t)v^{\frac{p-1}{p-1/2}} |\nabla v|^{s-2}
 \left[|A| + \frac{p -1}{2p-1} \frac{|\nabla v|^2}{v} \left( \frac{A (\xi, \xi)}{|A|}- \frac{{\rm tr} A}{|A|}\right) \right]^2\nonumber\\
&+ \frac{s p(p-1)^2}{(2p-1)^2} 
\left[ \frac{A (\xi, \xi)}{|A|}- \frac{{\rm tr} A}{|A|}\right]^2
\zeta(t) v^{\frac{2p}{1-2p}} |\nabla v|^{s+2}
 - \frac{sp \zeta(t)}{(2p-1)^2} v^{\frac{2p}{1-2p}} |\nabla v|^{s+2}\nonumber\\
&+ \zeta'(t) |\nabla v|^s - s\zeta(t) |\nabla v|^{s-2}\left[\frac{1}{2}\partial_t g 
+ pv^{\frac{p-1}{p-1/2}} {\mathscr Ric}_f^m(g)\right] ( \nabla v, \nabla v)\nonumber\\
&- \frac{sp}{2} \zeta(t) v^{\frac{p-1}{p-1/2}} \langle \nabla |\nabla v|^{s-2} , \nabla |\nabla v|^2 \rangle 
+ \zeta(t) s |\nabla v|^{s-2}\langle \nabla v, \nabla \Sigma^\star(t,x,v) \rangle\nonumber\\
&-sp\zeta(t)v^{\frac{p-1}{p-1/2}} |\nabla v|^{s-2} \left[\frac{\langle \nabla f, \nabla v \rangle}{\sqrt{m-n}}
+ \frac{2(p -1)}{2p-1} |\nabla v|^2 \frac{\sqrt{m-n}}{2v}\right]^2 
+ \Gamma'(v) \Sigma^\star(t, x, v) \nonumber\\
&+\frac{sp(p-1)^2(m-n)}{(2p-1)^2} \zeta(t)v^{\frac{2p}{1-2p}} |\nabla v|^{s+2} 
- p v^{\frac{1}{1-2p}} [\Gamma'(v)+v\Gamma''(v)]|\nabla v|^2. 
\end{align}
Now by ignoring non-positive terms 
and making use of the variational matrix inequality \eqref{matrix-variational-inequality} we conclude that 
\begin{align}
{\mathscr L}({\mathsf H}^{s}_\zeta[v])  - 
&\frac{2p^2v^{\frac{1}{1-2p}} }{2p-1} \langle \nabla v , \nabla {\mathsf H}^{s}_\zeta[v] \rangle
\le \frac{s p[(p-1)^2(m-1)-1]}{(2p-1)^2} \zeta(t) v^{\frac{2p}{1-2p}} |\nabla v|^{s+2} \nonumber\\
&+ \zeta'(t) |\nabla v|^s - s\zeta (t) |\nabla v|^{s-2}\left[\frac{1}{2}\partial_t g 
+ pv^{\frac{p-1}{p-1/2}} {\mathscr Ric}_f^m(g)\right] ( \nabla v, \nabla v)\nonumber\\
&-\frac{sp}{2} \zeta(t) v^{\frac{p-1}{p-1/2}} \langle \nabla |\nabla v|^{s-2} , \nabla |\nabla v|^2 \rangle 
+ s\zeta(t) |\nabla v|^{s-2}\langle \nabla v, \nabla \Sigma^\star(t,x,v) \rangle \nonumber\\
&+ \Gamma'(v) \Sigma^\star(t, x, v) - p v^{\frac{1}{1-2p}} [\Gamma'(v)+v\Gamma''(v)]|\nabla v|^2.
\end{align}

An application of inequality \eqref{eq11c-2.n} along with 
$\langle \nabla v, \nabla \Sigma^\star \rangle = \langle \nabla v, \Sigma^\star_x \rangle + \Sigma^\star_v |\nabla v|^2$ 
and $\langle \nabla |\nabla v|^{s-2}, \nabla |\nabla v|^2\rangle \ge 0$ when $s \ge 2$ 
gives the desired conclusion. 
\end{proof}

 \begin{corollary} \label{first-cor-last}
Let $(M,g,e^{-f} dv_g)$ be a smooth metric measure space with $M$ closed. 
Let $u$ be a positive smooth solution to \eqref{eq11} where $p_0(m)<p<1$ and set $v = u^{p-1/2}$. 
Assume the metric and potential satisfy the super flow \eqref {eq11c-2.n} with $\mathsf{k} \ge 0$. Moreover suppose 
the following conditions hold for some $a$: 
\begin{itemize}
\item $\Sigma^\star(t,x,v) \le a$, 
\item $\Gamma'(v) \Sigma^\star (t,x,v) \le 0$,
\item $\Gamma'(v) + v \Gamma''(v) \ge 0$,
\item $\langle \nabla v, \Sigma^\star_x(t,x,v) \rangle \le 0$.  
\end{itemize}
Then for all $x \in M$ and $0<t \le T$ we have 
\begin{equation}
|\nabla v|^s (x,t) 
\le e^{s(\mathsf{k}+a)t} \left\{ \max_{M} 
\left[ |\nabla v|^s + \Gamma(v) \right]_{t=0} - \Gamma(v(x,t)) \right\}. 
\end{equation}
\end{corollary}

\begin{proof}
Using \eqref{EQ-eq-4.4-2.n-2} and the assumptions on $\Gamma$ and $\Sigma$ in the statement of the corollary, we can write for 
${\mathsf H} = {\mathsf H}^s_\zeta[v]= \zeta(t)|\nabla v|^s + \Gamma(v)$
\begin{align}  \label{Eq-7.16}
{\mathscr L} & ({\mathsf H}) - \frac{2p^2v^{\frac{1}{1-2p}}}{2p-1} \langle \nabla v , \nabla {\mathsf H} \rangle\nonumber\\
\le &~[\zeta'(t)+ s \mathsf k \zeta(t) + s \zeta(t) \Sigma_v^\star(t,x,v)] |\nabla v|^s 
+sp \zeta(t) \frac{(p-1)^2(m-1)-1}{(2p-1)^2} v^{\frac{2p}{1-2p}} |\nabla v|^{s+2} \nonumber\\
&+ s \zeta(t) \langle \nabla v,  \Sigma^\star_x(t,x,v) \rangle |\nabla v|^{s-2}
+ \Gamma'(v) \Sigma^\star(t,x,v) -p v^{\frac{1}{1-2p}} [\Gamma'(v)+v\Gamma''(v)] |\nabla v|^2 \nonumber \\
\le &~[\zeta'(t)+ s (\mathsf k + a) \zeta(t)] |\nabla v|^s 
+\zeta(t) \frac{s p[(p-1)^2(m-1)-1]}{(2p-1)^2} v^{\frac{2p}{1-2p}} |\nabla v|^{s+2}.
\end{align}
The function $\zeta(t) = e^{-s(\mathsf{k}+a)t}$ is non-negative, smooth and satisfies $\zeta'+s(\mathsf{k}+a)\zeta=0$. 
Thus substituting in \eqref{Eq-7.16} and noting that $(p-1)^2(m-1) \le 1$ as a result of $p_0<p<1$ we conclude that 
${\mathscr L}({\mathsf H}) - 2p^2 v^{1/(1-2p)}/(2p-1)  \langle \nabla v, \nabla {\mathsf H} \rangle \le 0$. 
The assertion is now a consequence of the weak maximum principle.
\end{proof}

\begin{corollary} \label{second-cor-last}
Let $(M,g,e^{-f} dv_g)$ be a smooth metric measure space with $M$ closed. 
Let $u$ be a positive smooth solution to \eqref{eq11} where $p_0(m)<p <1$ and set $v =u^{p-1/2}$. 
Assume the metric and potential satisfy the super flow \eqref {eq11c-2.n} with $\mathsf{k} \ge 0$. Moreover suppose 
the following conditions hold: 
\begin{itemize}
\item $\Sigma^\star(t,x,v) \le 0$,
\item $\Sigma^\star_v(t,x,v) \le 0$,
\item $\langle \nabla v, \Sigma^\star_x(t,x,v) \rangle \le 0$. 
\end{itemize}
Then for $x \in M$ and $0<t \le T$ we have 
\begin{equation}
\frac{t |\nabla v(x,t)|^2}{1+2{\mathsf k}t} \le \frac{(2p-1)^2}{8p^3} \left[ \max_M v^\frac{2p}{2p-1}(x,0) - v^\frac{2p}{2p-1}(x,t) \right]. 
\end{equation}
\end{corollary}

\begin{proof}
For ${\mathsf H}[v] = \zeta(t)|\nabla v|^2 + c v^\frac{2p}{2p-1}$ where 
$\Gamma(v)=cv^\frac{2p}{2p-1}$, we have from 
Lemma \ref{evolution Rpq-2.n-inequality-2},
\begin{align} \label{WMXprclelast}
{\mathscr L} & ({\mathsf H}) - \frac{2p^2v^{\frac{1}{1-2p}}}{2p-1} \langle \nabla v , \nabla {\mathsf H} \rangle\nonumber\\
\le &~[\zeta'(t)+ 2 \mathsf k \zeta(t) + 2 \zeta(t) \Sigma_v^\star(t,x,v)] |\nabla v|^2 
+2p \zeta(t) \frac{(p-1)^2(m-1)-1}{(2p-1)^2} v^{\frac{2p}{1-2p}} |\nabla v|^4 \nonumber\\
&+ 2\zeta(t) \langle \nabla v,  \Sigma^\star_x(t,x,v) \rangle 
+ \frac{2pc}{2p-1} v^\frac{1}{2p-1} \Sigma^\star(t,x,v) - \frac{8p^3c}{(2p-1)^2} |\nabla v|^2. 
\end{align}
Recalling the assumptions in the corollary and using $c=(2p-1)^2/(8p^3)$ then gives 
\begin{align} 
{\rm RHS \,} \eqref{WMXprclelast}
\le
[\zeta'(t) + 2\mathsf{k} \zeta(t) -1] |\nabla v|^2 
+ \frac{(p-1)^2(m-1)-1}{2p-1} \zeta(t) v^\frac{2p}{1-2p} |\nabla v|^4.  \nonumber 
\end{align}

Now taking $\zeta(t) = t/(1+2\mathsf{k} t)$ gives $(\zeta' + 2 \mathsf{k} \zeta -1) \le 0$. Therefore since $p_0<p<1$ 
we have ${\mathscr L}({\mathsf H}) - [2p^2v^{1/(1-2p)}/(2p-1)] \langle \nabla v, \nabla {\mathsf H} \rangle \le 0$. 
The conclusion now follows by an application of the weak maximum principle.  
\end{proof}

\begin{remark}
In both corollaries above we can admit the left end-point of the $p$ interval, i.e., $p=p_0(m)$ provided 
that the strict inequality $p>1/2$ is observed. This specifically means that for $m>5$ where 
$p_0(m) = 1-1/\sqrt{m-1}>1/2$ we can extend the range to $p_0(m) \le p <1$.
\end{remark}

\qquad \\
{\bf Acknowledgement.} The authors gratefully acknowledge support from the Engineering and Physical 
Sciences Research Council (EPSRC) through the grant EP/V027115/1.

\qquad \\
{\bf Data Availability.}  Data sharing is not applicable to this article as no datasets were generated or analysed during this study.


\begin{thebibliography}{99}

\bibitem{AM} E.~Acerbi, R.~Mingione, {\it Gradient estimates for a class of parabolic systems}, Duke Math. J., 136, (2007), 285-320.

\bibitem{AGS} L.~Ambrosio, N.~Gigli, G.~Savar\'e, Gradient Flows: In Metric Spaces and the Space of Probability Measures, 
Lectures in Mathematics, ETH Z\"urich, 2nd Edition, Birkh\"auser, 2008.

\bibitem{Aub} T.~Aubin, {\it Nonlinear Analysis on Manifolds}, Springer, New York 1982.

\bibitem{AronBen} D.G.~Aronson, P.~B\'enilan, {\it R\'egularit\'e des solutions de l'\'equation des milieux poreux dans 
${\mathbb R}^n$}, C.R. Acad. Sci. Paris S\'er. A-B, 288, (1979), 103--105. 

\bibitem{Aron} D.G.~Aronson, {\it The porous medium equation}, In: Some Problems in Nonlinear Diffusion, 
A.~Fasano and M.~Primicerio, Eds., Lecture Notes in Mathematics, {\bf 1224}, Springer-Verlag, (1986), 1--46.

\bibitem {BaCP}{M.~B\v{a}ile\c{s}teanu, X. Cao, A. Pulemotov},  
{\it Gradient estimates for the heat equation under the Ricci flow}, J. Funct. Anal., 258 (2010), 3517--3542. 

\bibitem{Bak}{D.~Bakry, I.~Gentil, M.~Ledoux}, Analysis and Geometry of Markov Diffusion Operators, 
A Series of Comprehensive Studies in Mathematics, {\bf 348}, Springer, 2012.


\bibitem{BBDGV} {A.~Blanchet, M.~Bonforte, J.~Dolbeault, G.~Grillo, J.L.~V\'azquez}, 
{\it Asymptotics of the fast diffusion equation via entropy estimates}, Arch. Rational Mech. Anal., 191 (2009), 347--385. 

\bibitem{BDM} V.~B\"ogelein, F.~Duzaar, G.~Mingione, {\it The regularity of general parabolic systems with degenerate diffusion}, 
Mem. Amer. Math. Soc., 221 (1041), (2013), vi+143~pp., AMS.  

\bibitem{BonBon} M.~Bonforte, J.~Dolbeault, G.~Grillo, J.L.~V\'azquez, {\it Sharp rates of decay of solutions to the nonlinear 
fast diffusion equation via functional inequalities}, PNAS, 107(38) (2010), 16549--16464.

\bibitem{Bon-Dol} M.~Bonforte, J.~Dolbeault, B.~Nazaret, N.~Simonov, {\it Stability in Gagliardo-Nirenberg-Sobolev Inequalities:
Flows, Regularity and the Entropy Method}, Mem. Amer. Math. Soc., To appear. 

\bibitem{CCL} E.A.~Carlen, J.A.~Carrillo, M.~Loss, {\it Hardy-Littlewood-Sobolev inequalities via fast diffusion}, 
PNAS, 107(46) (2010), 19696--19701.  


\bibitem{Chow}{B.~Chow, P.~Lu, L.~Nei}, Hamilton's Ricci Flow, Graduate Studies in Mathematics {\bf 77}, AMS, 2006.	

\bibitem{DaskKe}{P.~Daskalopoulos, C.~Kenig}, Degenerate Diffusions: Initial Value Problems and Local Regularity Theory, 
Tracts in Mathematics {\bf 1}, European Mathematical Society, 2007.

\bibitem{Giaq} M.~Giaquinta, Multiple Integrals in the Calculus of Variations and Nonlinear Elliptic Systems, 
Annals of Mathematics Studies, Vol.~{\bf 105}, Princeton University Press, 1983. 
 
\bibitem{Ha93}{R.~Hamilton},  {\it A matrix Harnack estimate for heat equation}, Comm. Anal. Geom., (1993), 113--126.	 
 
\bibitem{JiCh}{X.~Jiang, Y.~Cheng}, {\it Gradient estimates for fast diffusion equations on Riemannian manifolds}, 
J.~Math.~Anal.~Appl., 472 (2019), 1369--1376. 
 
\bibitem{Hunag-Li} G.~Huang, H.~Li, {\it Gradient estimates and entropy formulae of porous medium and fast diffusion 
equations for the Witten Laplacian}, Pac.~J.~Math., 268 (2014), 47--78. 

\bibitem{JKO} R.~Jordan, D.~Kinderleher, F.~Otto, {\it The variational formulation of the Fokker-Planck equation}, 
SIAM J. Math. Anal., 29 (1), (1998), 1--17.

\bibitem{KT}{J.~Kristensen, A.~Taheri}, {\it Partial regularity of strong local minimizers in the calculus of variations}, 
Arch. Rational Mech. Anal., 170, (2003), 63--89.

\bibitem{[LiP]}{P.~Li}, Geometric Analysis, Cambridge Studies in Advanced Mathematics, {\bf 134}, CUP, 2012.

\bibitem{[LY86]}{P.~Li, S.T.~Yau}, {\it On the parabolic kernel of Schr\"odinger operator}, Acta Math., 156 (1986), 153--201.
 
\bibitem{Lott}{J.~Lott}, {\it Some geometric properties of the Bakry-\'Emery Ricci tensor}, Comment. Math. Helv., 78 (2003), 865--883. 

\bibitem{PLu} P.~Lu, L.~Ni, J.~V\'azquez, C.~Villani, {\it Local Aronson-B\'enilan estimates and entropy formulae for porous 
medium and fast diffusion equations on manifolds}, J.~Math.~Pures~Appl., 91 (2009), 1--19.

\bibitem{MZS}{L.~Ma, L.~Zhao, X.~Song}, {\it Gradient estimates for the degenerate parabolic equation 
$u_t=\Delta F(u) + H(u)$ on manifolds}, J.~Diff. Eq., 244 (2008), 1157-1177.

\bibitem{Otto} F.~Otto, {\it The geometry of dissipative evolution equations: The porous medium equation}, 
Comm. Partial Differential Equations, 26 (2001), No.~1-2, 101--174. 

\bibitem{Pe02}{G.~Perelman}, {\it The entropy formula for the Ricci Flow and its geometric application}, arXiv: math. DG/0211159v1 (2002).

\bibitem{SZ} {P.~Souplet, Q.S.~Zhang} {\it Sharp gradient estimate and Yau's Liouville theorem for the heat equation on noncompact manifolds}, 
Bull. Lond. Math. Soc., 38 (2006), 1045--1053. 

\bibitem{Taheri-book-one} A.~Taheri, 
Function Spaces and Partial Differential Equations, Vol.~{\bf I}, Oxford Lecture Series in Mathematics and its Applications, {\bf 40}, OUP, 2015.   

\bibitem{Taheri-book-two} A.~Taheri, 
Function Spaces and Partial Differential Equations, Vol.~{\bf II}, Oxford Lecture Series in Mathematics and its Applications, {\bf 41}, OUP, 2015. 

\bibitem{Taheri-GE-1} A.~Taheri, {\it Liouville theorems and elliptic gradient estimates for a nonlinear parabolic equation involving the Witten Laplacian}, 
Published online in: Adv. Calc. Var., De Gruyter, 2021.

\bibitem{Taheri-GE-2} A.~Taheri, {\it Gradient estimates for a weighted $\Gamma$-nonlinear parabolic equation coupled with a super Perelman-Ricci flow 
and implications}, Published online in: Potential Anal., Springer, 2021. 

\bibitem{TVahNA} A.~Taheri, V.~Vahidifar, {\it Gradient estimates for a nonlinear parabolic equation on smooth metric measure spaces with 
evolving metrics and potentials}, Nonlinear Anal., 232, 2023. 

\bibitem{TVahCurv} A.~Taheri, V.~Vahidifar, {\it Curvature conditions, Liouville-type theorems and Harnack inequalities for a nonlinear 
parabolic equation on smooth metric measure spaces}, Adv. Nonlinear Studies, De Gruyter, 2024.

\bibitem{TVDiffHar} A.~Taheri, V.~Vahidifar, {\it Differential Harnack estimates for a weighted nonlinear parabolic equation under a super 
Perelman-Ricci flow and implications}, To appear in: Proc. Roy. Soc. Edin., Mathematics A, CUP, 2024.

\bibitem{TVah-PME} A.~Taheri, V.~Vahidifar, {\it The nonlinear porous medium equation on smooth metric measure spaces and the 
Ricci-Perelman super flow}, Submitted for publication 2024. 

\bibitem{Vaz} J.~V\'azquez, The Porous Medium Equation, Oxford Mathematical Monographs, Clarendon Press, OUP, 2007.

\bibitem{VC} {C.~Villani}, Optimal transport: Old and New, A Series of Comprehensive Studies in Mathematics, {\bf 338}, Springer, 2008.

\bibitem{Wang}{F.Y.~Wang}, {\it Analysis for Diffusion Processes on Riemannian Manifolds}, Advanced Series on Statistical Science $\&$ Probability, 
Vol.~{\bf 18}, World Scientific, 2013. 

\bibitem{XXu}{X.~Xu}, {\it Gradient estimates for $u_t=\Delta F(u)$ on manifolds and some Liouville-type theorems}, J.~Diff.~Eq, 252 (2012), 1403-1420.

\bibitem{Zhang} Q.S.~Zhang, Sobolev inequalities, heat kernels under Ricci flow and the Poincar\'e conjecture, CRC Press, 2011. 

\bibitem{Zhu}{X.~Zhu}, {\it Hamilton's gradient estimates and Liouville theorems for fast diffusion equations on noncompact Riemannian manifolds}, 
Proc. Amer. Math. Soc., 139 (2011), 1637--1644.

\end{thebibliography}
\end{document}